\documentclass{article}

\usepackage{verbatim,amsfonts,amsthm,epsfig,fullpage,amsmath,psfig,color}

\begin{document}

\newcommand{\chek}[1]{{{#1}}}
\newcommand{\dima}[1]{{{#1}}}

\newtheorem{definition}{Definition}
\newtheorem{lemma}{Lemma}
\newtheorem{claim}{Claim}
\newtheorem{theorem}{Theorem}
\newtheorem{corollary}{Corollary}
\newtheorem{conjecture}{Conjecture}
\newtheorem{fact}{Fact}
\newtheorem{proposition}[lemma]{Proposition}

\def\whp{w.h.p.}
\def\wupp{with uniformly positive probability}
\def\Whp{W.h.p.}
\def\eg{e.g.\ }
\def\ie{i.e.\ }

\newcommand{\p}{p}
\newcommand{\q}{q}
\renewcommand{\log}{\ln}

\newcommand{\viol}{\not\models}
\newcommand{\ra}{\rightarrow}

\newcommand{\e}{\epsilon}
\newcommand{\dd}{\delta}

\newcommand{\rc}{{c}}
\newcommand{\fc}{{\omega}}

\newcommand{\Z}{Z}

\newcommand{\up}{*}
\newcommand{\rr}{T}
\newcommand{\ttk}{t_k}

\newcommand{\amax}{\alpha_{\max}}
\newcommand{\gmax}{g_{\max}}
\def\a{\alpha}

\def\e{\varepsilon}
\def\FF{F}

\def\ex{{\mathbf E}}
\def\E{{\mathbf E}}
\newcommand{\one}{{\mathbf{1}}}
\newcommand{\PS}{\mathbf{P}}
\newcommand{\Pa}{{\mathbf P_{\sigma}}}
\newcommand{\Ea}{{\mathbf E_{\sigma}}}

\title{On the Maximum Satisfiability of Random Formulas}

\author{
    Dimitris Achlioptas\thanks{Part of this work was done while visiting UC Berkeley.} \\
    Microsoft Research, Redmond, Washington \\
    {\tt optas@microsoft.com}
    \and
    Assaf Naor \\
    Microsoft Research, Redmond, Washington \\
    {\tt anaor@microsoft.com}
    \and Yuval Peres\thanks{Research supported by NSF Grant
DMS-0104073 and a Miller Professorship at UC Berkeley.}\\
Departments of Statistics and Mathematics, University of California, Berkeley\\
{\tt peres@stat.berkeley.edu} }

\date{\empty}

\maketitle

\begin{abstract}
Maximum satisfiability is a canonical NP-hard optimization problem
that \dima{appears empirically hard for random instances}. In
particular, its apparent hardness on random $k$-CNF formulas
\dima{of certain densities} was recently suggested by Feige as a
starting point for studying inapproximability. At the same time,
it is rapidly becoming a canonical problem for statistical
physics. In both of these realms, evaluating new ideas relies
crucially on knowing the maximum number of clauses one can
typically satisfy in a random $k$-CNF formula. In this paper we
give asymptotically tight estimates for this quantity.
Specifically, let us say that a $k$-CNF  is {\em $\p$-satisfiable}
if there exists a truth assignment satisfying $1-2^{-k}+\p 2^{-k}$
of all clauses (observe that every $k$-CNF is 0-satisfiable).
Also, let $F_k(n,m)$ denote a random \mbox{$k$-CNF} on $n$
variables formed by selecting uniformly and independently $m$ out
of all $2^k \binom{n}{k}$ possible $k$-clauses.

Let \dima{$\tau(p) = 2^k \ln 2/({\p+(1-\p)\log (1-\p)})$}. It is
easy to prove that for every $k \geq 2$ and every $\p \in (0,1]$,
if $r \geq \tau(\p)$ then the probability that $F_k(n,rn)$ is
$\p$-satisfiable tends to 0 as $n \to \infty$. We prove that
\dima{there exists a sequence $\delta_k \to 0$ such that if $r
\leq (1-\delta_k) \tau(\p)$ then the probability that $F_k(n,rn)$
is $\p$-satisfiable tends to 1 as $n \to \infty$. The sequence
$\delta_k$ tends to 0 exponentially fast in $k$. Indeed, even for
moderate values of $k$, \eg $k=10$, our result gives very tight
bounds for the number of satisfiable clauses in a random $k$-CNF.
In particular, for $k>2$ it improves upon all previously known
such bounds.}
\end{abstract}

\newpage

\section{Introduction}
Given a Boolean CNF formula $F$, the Satisfiability problem is to
determine whether there exists a truth assignment that satisfies
$F$. When $F$ has exactly $k$ literals in each clause,
Satisfiability is known as $k$-SAT and is NP-complete~\cite{cook}
for all $k \geq 3$. A natural generalization of satisfiability is
determining whether there exists a truth assignment that satisfies
a given number of clauses in $F$. For $k$-CNF this problem is
known as Max $k$-SAT and is NP-complete for all $k\geq 2$
(see~\cite{GaJo79}).

Optimization problems with random inputs are pervasive in
operations research (e.g., the travelling salesman problem and
variants), in statistical physics (determining ground states of
spin glasses) and in computer science. An interesting source of
Max $k$-SAT instances comes from considering $k$-CNF chosen
uniformly at random (see below). Historically, the motivation for
studying such formulas has been the desire to understand the
hardness of ``typical'' instances.
 Random $k$-CNF are
by now the most studied generative model for random formulas and
have been a very popular benchmark for testing and tuning
satisfiability algorithms. In fact, some of the better practical
ideas in use today come from insights gained by studying the
performance of algorithms on random
$k$-CNF~\cite{slm:gsat,sk:walksat,heavytail}.

A natural starting point for considering Max $k$-SAT is the
observation that for every $k$-CNF formula there exists a truth
assignment satisfying at least $(1-2^{-k})$ of all clauses.
Indeed, if such a formula has $m$ clauses, the average over all
$2^n$ truth assignments of the number of satisfied clauses is
precisely $(1-2^{-k})m$. With this in mind, we will say that a
$k$-CNF formula is {\em $\p$-satisfiable}, where $p \in [0,1]$, if
there exists a truth assignment satisfying $1-2^{-k}+\p 2^{-k}$ of
all clauses.

\dima{To consider random $k$-CNF formulas, let $C_k$ denote the
set of all $(2n)^k$ possible disjunctions of $k$ literals on some
canonical set of $n$ Boolean variables. To form a random $k$-CNF
formula $F_k(n,m)$ with $m$ clauses we select uniformly,
independently and with replacement $m$ clauses from $C_k$ and take
their conjunction\footnotemark[2].} We will say that a sequence of
random events ${\mathcal E}_n$ occurs {\em with high probability\/
} (w.h.p.) if $\lim_{n \ra \infty} \Pr[{\mathcal E}_n]=1$ and {\em
with uniformly positive probability\/ } if $\liminf_{n \ra \infty}
\Pr[{\mathcal E}_n]>0$. We emphasize that throughout the paper $k$
is arbitrarily large but fixed, while $n \rightarrow \infty$. For
every $k \geq 2$ and $\p \in (0,1]$, let
\begin{eqnarray*}
r_k(p)  & \equiv &
        \sup\{r : F_k(n,rn)\mbox{ is $p$-satisfiable    \whp}\} \\
        & \leq &
        \inf\{r : F_k(n,rn)\mbox{ is {\em not\/} $p$-satisfiable  \whp}\}
        \; \equiv \; r_k^\up(p) \enspace .
\end{eqnarray*}

\footnotetext[2]{Our discussion and results hold in all common
models for random $k$-CNF, \eg when clause replacement is not
allowed and/or when each $k$-clause is formed by selecting $k$
distinct, non-complementary literals with/without ordering. The
model defined here is best suited for our calculations. We further
comment on its relationship to other models in the end of
Section~\ref{outline}.}

One of the most intriguing aspects of random formulas is the {\em
Satisfiability Threshold Conjecture} which asserts that
$r_k(1)=r_k^\up(1)$ for every $k \geq 3$.  Much work has been done
to bound $r_k(1)$ and $r_k^\up(1)$. Currently, the best rigorous
bounds for general $k \geq 3$, from~\cite{yuval,DuBo}
respectively, are: $2^k \ln 2 - O(k) < r_k \le r_k^\up< 2^k \ln 2
- O(1)$. \dima{For $p<1$, the bounds for $r_k(\p), r_k^\up(\p)$
were much further apart.}

The state of the art for general $k$ was presented in an important
recent paper by Coppersmith, Gamarnik, Hajiaghayi, and
Sorkin~\cite{sork}, where it was proved \dima{(see~\eqref{sor_or}
for a more precise formulation)} that there exists an absolute
constant $c>0$ such that for all $k$ and all $p \in (0, p_0(k)]$,
\begin{equation} \label{sor2}
\frac{c}{k} \,  \frac{2^{k+1} \ln 2}{\p^{2}} \le r_k(\p) \le
r_k^\up(\p) \leq \frac{2^{k+1} \ln 2}{\p^{2}(1+o(1))} \enspace .
\end{equation}


The upper bound in~\eqref{sor2} was proved via the first moment
method, while the lower bound is algorithmic. For small $k$ the
two are reasonably close, but the ratio between them tends to
infinity as $k$ grows; this naturally raises the question which
bound is closer to the truth. 
Our main result \dima{resolves this question by pinpointing the
values of $r_k(\p)$ and $r_k^\up(\p)$ with relative error that
tends to zero exponentially fast in $k$.} For every $\p \in (0,1)$
denote
\begin{equation}\label{defrr}
\rr_k(\p)=\frac{2^k\ln 2}{\p+(1-\p)\log (1-\p)} \,,
\end{equation}
and let $\rr_k(1)=2^k\ln 2$ so that $\rr_k(\cdot)$ is continuous
on $(0,1]$.
\begin{theorem}\label{main}
There exists a  sequence $\delta_k=O(k 2^{-k/2})$, such that for
all $k \geq 2$ and $p \in (0,1]$,
\begin{equation}\label{eq:main}
 (1-\delta_k) \, {\rr}_k(\p) < r_k(\p)\le r_k^{\up}(\p)\le {\rr}_k(\p) \enspace .
\end{equation}
\end{theorem}
\dima{The upper bound in (\ref{eq:main}) follows from well-known
tail estimates}. Taylor expansion gives that as $p \to 0$,
$$
\rr_k(\p)=\frac{2^k\ln 2}{\p^2/2+O(\p^3)} \enspace ,
$$
so as $p \to 0$, we can sharpen \eqref{sor2} to
\begin{equation} \label{nsor2}
 (1-\delta_k) \frac{2^{k+1} \ln 2}{\p^{2}+O(p^3)} \le r_k(\p) \le
r_k^\up(\p) \leq \frac{2^{k+1} \ln 2}{\p^{2}+O(p^3)} \enspace .
\end{equation}

Our proof  of Theorem~\ref{main} actually yields an explicit lower
bound for $r_k(p)$ for each $k \geq 2$.
For $k =2$, \ie Max 2-SAT, the algorithm presented in~\cite{sork}
dominates our lower bound uniformly, \ie for every density it
satisfies a greater fraction of all clauses. Already for $k \geq
3$, though, our methods yield a better bound, as indicated by the
following plots.

\vspace*{0.7cm}

\begin{minipage}{3in}
\begin{center}
        \centerline{\hbox{
        \psfig{figure=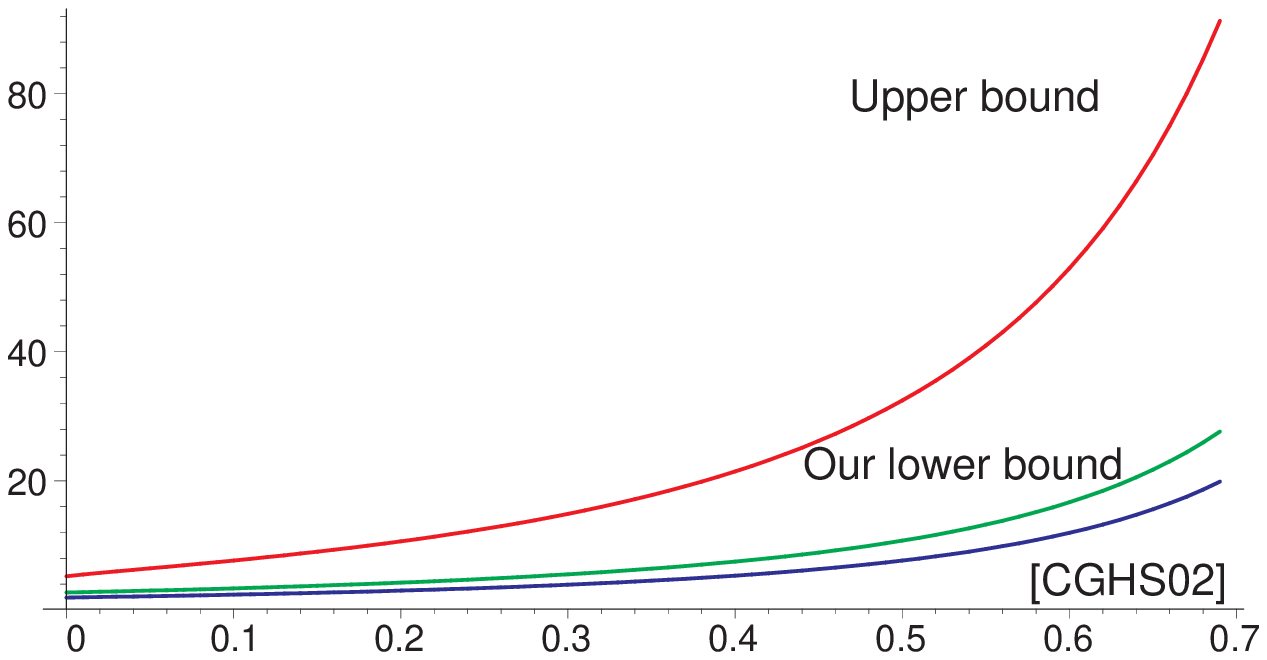,height=1.5in,width=2in}
        }}
        \centerline{$k=3$}
\end{center}
\end{minipage}\    \
\begin{minipage}{3in}
\begin{center}
        \centerline{\hbox{
        \psfig{figure=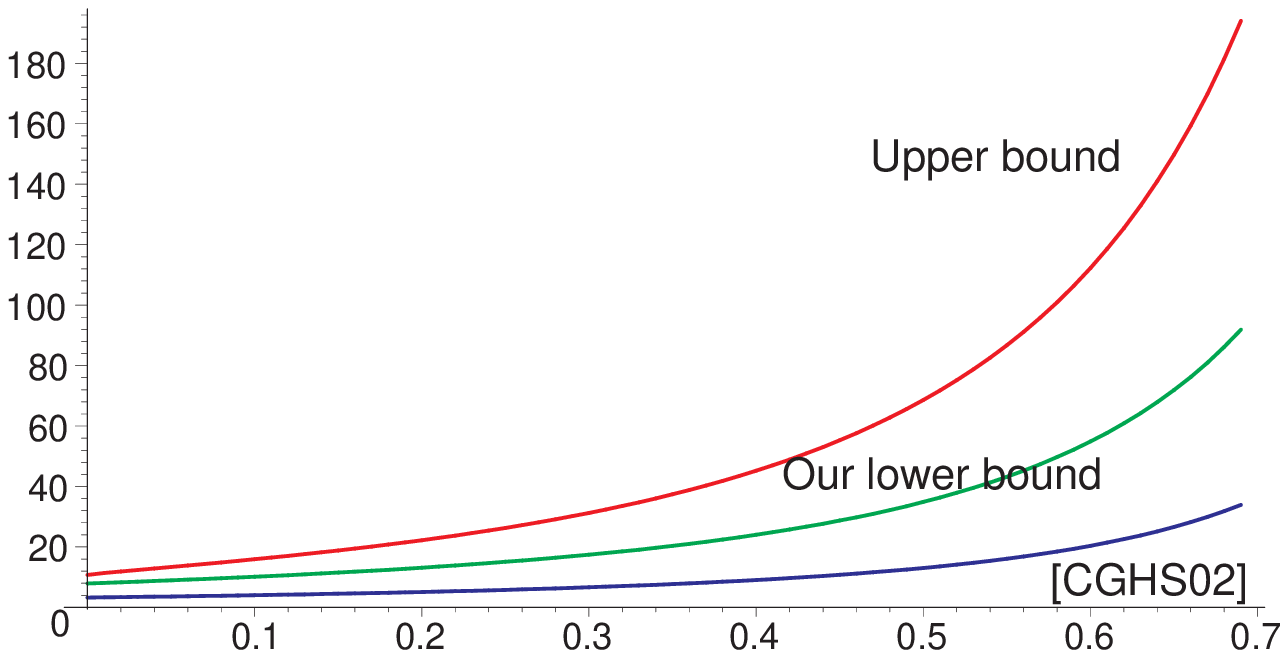,height=1.5in,width=2in}
        }}
        \centerline{$k=4$}
\end{center}
\end{minipage}

\vspace*{0.6cm}

\begin{minipage}{3in}
\begin{center}
        \centerline{\hbox{
        \psfig{figure=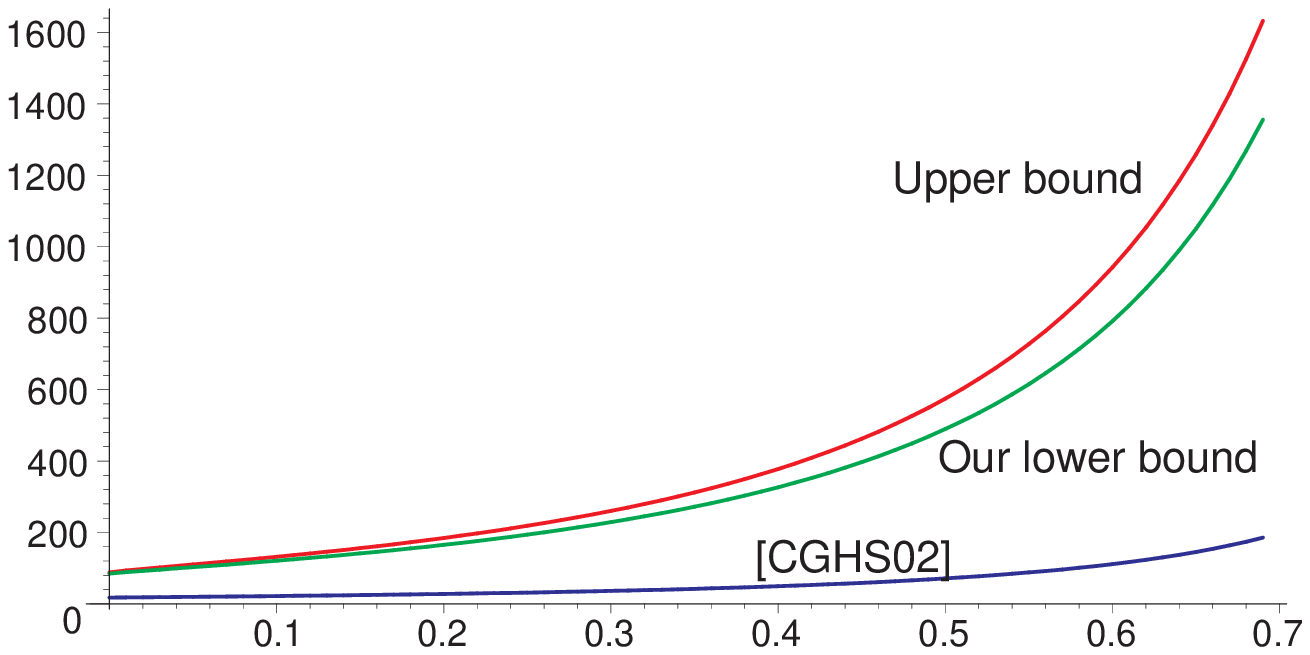,height=1.5in,width=2in}
        }}
        \centerline{$k=7$}
\end{center}
\end{minipage}\    \
\begin{minipage}{3in}
\begin{center}
        \centerline{\hbox{
        \psfig{figure=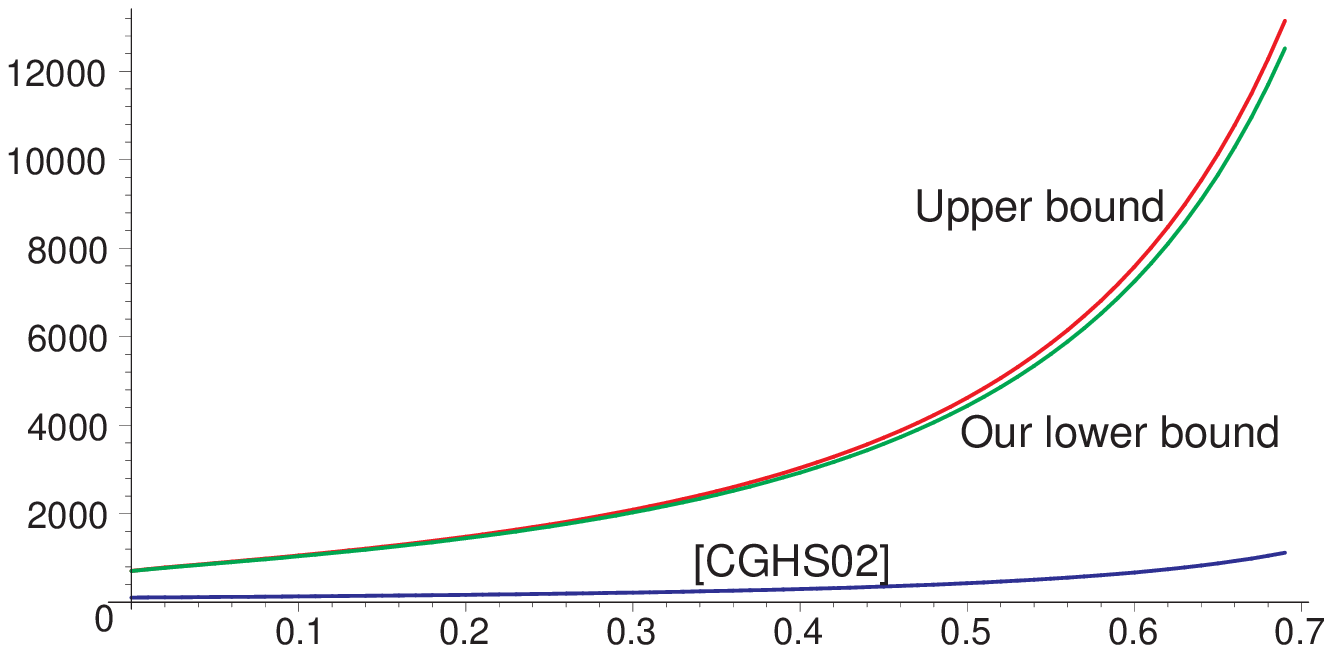,height=1.5in,width=2in}
        }}
        \centerline{$k=10$}
\end{center}
\end{minipage}

\vspace*{-0.2cm}

\begin{center}
Figure 1. Upper and lower bounds for the density $r$ as a function
of $\q=1-p$.
\end{center}

Our approach in proving Theorem~\ref{main} is non-algorithmic,
based instead on a delicate application of the second moment
method to a random generating function in two variables. It is
notoriously difficult to obtain precise asymptotics from such
random multivariable generating functions; the fact that this is
possible for random Max $k$-SAT is technically due to the
surprising cancellation of four terms of equal magnitude in our
analysis, leaving only lower order terms. This cancellation hints
at the existence of some unexpected hidden structure in random Max
$k$-SAT; characterizing this structure combinatorially (rather
than just analytically) appears to us  worthy of further study.

\subsection{Background}

For a random formula $F_k(n,m)$, denote by $s_k(n,m)$ the random
variable equal to the maximum (over all truth assignments
$\sigma$) of the number of clauses satisfied by $\sigma$. Perhaps
the first rigorous study of random Max $k$-SAT appeared in the
work of Frieze, Broder and Upfal~\cite{BFU} where  it was shown
that $s_k(n,m)$ is sharply concentrated
around its mean. Specifically, 
\begin{theorem}[{\cite{BFU}}]\label{conc}
$\displaystyle{ \Pr\biggl[\bigl|s_k(n,m) - \ex[s_k(n,m)]\bigr| > t
\biggr] < 2 \exp(-2t^2/m)}$.
\end{theorem}

The following corollary allows us to infer high probability
results from positive probability results.
\begin{corollary}\label{booster}
If $F_k(n,rn)$ is $\p_0$-satisfiable \wupp, then for every
constant $\p<\p_0$, $F_k(n,rn)$ is $\p$-satisfiable \whp\
\end{corollary}
\begin{proof}
Let $S \equiv (1-2^{-k}+\p_0 2^{-k})rn$. Since $F_k(n,rn)$ is
$\p_0$-satisfiable \wupp, $\ex[s_k(n,rn)] > S - n^{2/3}$. For,
otherwise, Theorem~\ref{conc} would imply that the probability of
$\p_0$-satisfiability is exponentially small. By the same token,
$\Pr[s_k(n,rn) < S- 2  n^{2/3}] =o(1)$, implying the claim.
%
%
\end{proof}

Thus it will suffice to find, for every $\p \in (0,1]$, a value
$r=r(\p)$ such that $F_k(n,rn)$ is $\p$-satisfiable with uniformly
positive probability and rely on Corollary~\ref{booster} to get a
high probability result.

\dima{Regarding the mean in Theorem~\ref{conc},} in view of the a
priori bound $s_k(n,m) \ge (1-2^{-k})m$, it is natural to consider
$\Phi_k(n,m)=\E[s_k(n,m)]-(1-2^{-k})m$, measuring how much the
optimum truth assignment does better than the a priori bound in
expectation (over random $k$-CNF). In~\cite{sork} it was shown
that for all $k$, for sufficiently large $r$, as $n \to \infty$
one has in $F_k(n,rn)$
\begin{equation} \label{sor1}
\frac{2}{k+1}\sqrt{\frac{k}{\pi 2^k}}\times {\sqrt{r}}-O(1) \le
\frac{\Phi_k(n,rn)}{n} \le \sqrt{\frac{(2^k-1) \ln 2}{2^{2k-1}}}
\times \sqrt{r}  \enspace .
\end{equation}
This is equivalent to the assertion that for $\p$ sufficiently
small,
\begin{equation}\label{sor_or}
\frac{k2^{k+2}}{\pi (k+1)^2}\times p^{-2} -O(\p^{-1}) \le r_k(\p)
\le r_k^\up(\p) \le 2(2^k-1) \ln 2
 \times \p^{-2} \enspace ,
\end{equation}
\dima{which is a more precise formulation of~\eqref{sor2}.}

Since for $k=2$, the threshold for satisfiability is known, namely
$r_2(1)=r_2^{\up}(1)$=1, in~\cite{sork} very fine results were
derived for $s_k(n,rn)$ when $r \approx 1$. In particular, when
$r=1+\varepsilon$ one has $\E[s_2(n,m)] =
(1+\varepsilon-O(\varepsilon^3))n$, while for large $r>1$ the
bound in~\eqref{sor1} can be improved to
\[
\frac{\sqrt{8}-1}{3\sqrt{\pi}} \times {\sqrt{r}} -O(1) \le
\frac{\Phi_2(n,rn)}{n} \le \sqrt{\frac{3 \ln 2}{8}} \times
{\sqrt{r}} \enspace .
\]

\dima{Another intriguing aspect of random $k$-CNF formulas is
their proof complexity. In a seminal paper, Chv\'{a}tal and
Szemer\'{e}di~\cite{ChvSz} proved that for all $k \geq 3$ and
$r>2^k \ln 2$ there exists $\e=\e(r)$ such that \whp\ every
resolution refutation of $F_k(n,rn)$ contains at least $(1+\e)^n$
clauses. Since then there have been a number of extensions of this
result~\cite{bp:resolve,BKPS} and it is widely believed that
random $k$-CNF are hard for much stronger proof systems than
resolution. Indeed, recently, Feige~\cite{feiger} showed that a
hypothesis asserting that proving unsatisfiability of random
$k$-CNF with $r \gg 2^k \ln 2$ is hard, implies a number of strong
inapproximability results. A closely related hypothesis is that
approximating Max $k$-SAT for such formulas is also hard for all
$k \geq 2$. Recent work by Fernandez De la Vega and
Karpinski~\cite{karpinski} proves that one can approximate Max
3-SAT on $F_3(n,rn)$ within $9/8$ which is better than the trivial
$8/7$ bound.}

\section{Outline}\label{outline}

\subsection{Understanding correlation sources in MAX $k$-SAT}
The following easy consequence of the Cauchy-Schwarz inequality
underlies the second moment method.
\begin{lemma}\label{lemma:sec}
For any non-negative random variable $X$,
\begin{equation}\label{eq:smm}
\Pr[X > 0] \,\ge\, \frac{\ex[X]^2}{\ex[X^2]} \enspace .
\end{equation}
\end{lemma}

Thus, for any fixed $\p \in (0,1]$ one can let $X$ denote the
number of $\p$-satisfying assignments and
apply~\eqref{lemma:sec}  to bound \mbox{$\Pr[X>0]$} from below.
Unfortunately, it turns out that for every
$r>0$, there exists a constant $\beta=\beta(k,r)>1$ such that $\ex[X^2]
> \beta^n \ex[X]^2$. As a result, this straightforward approach
only gives a trivial lower bound on the probability of
$\p$-satisfiability.

In~\cite{yuval}, it was shown that in the case $\p=1$ a major
factor in the excessive correlations behind the above failure is
the following form of populism: leaning toward the majority vote
truth assignment. To see this, first observe that truth
assignments that satisfy more literal occurrences than average,
have higher probability of being satisfying. At the same time, in
order to satisfy many literal occurrences such assignments tend to
agree with each other (and the majority truth assignment) on more
than half the variables. As a result, the successes of such
assignments tend to be highly correlated, thus dominating
$\ex[X^2]$. In order to avoid this pitfall, we would like, as in
\cite{yuval}, to apply the second moment method to truth
assignments that satisfy, approximately, half of all literal
occurrences; we call such truth assignments  {\em ``balanced''}.
In the context of $\p$-satisfiability, however, there are new
obstacles \chek{to overcome before obtaining} a lower bound for
$r_k(p)$ that asymptotically matches the upper bound.
To capture the behavior of balanced truth assignments we begin by
defining two ``fitness'' gauges.\medskip

Given any $k$-CNF formula ${\FF}$ on $n$ variables
and any truth assignment $\sigma \in \{0,1\}^n$ let
\begin{enumerate}
\item
$H= H(\sigma,{\FF})$ be the number of satisfied {literal
occurrences} in ${\FF}$ under $\sigma$, minus the number of
unsatisfied literal occurrences in ${\FF}$ under $\sigma$.
\item
$U= U(\sigma,{\FF})$ be the number of unsatisfied clauses in
${\FF}$ under $\sigma$.
\end{enumerate}

We would like to focus on truth assignments that are
balanced and $p$-satisfying, up to fluctuations
one would expect from a central limit theorem, i.e.,  truth assignments
$\sigma$ such that
\begin{align}
|H(\sigma,F)| & \leq A\sqrt{m} \label{iopu}\\
|U(\sigma,{\FF}) (1-\p)2^{-k} m| & \leq A\sqrt{m} \label{colopu} \enspace .
\end{align}

\chek{To do this let us write} $u_0 \equiv (1-\p)2^{-k}$ and fix
some $\gamma,\eta < 1$. \chek{Now,} for a random $k$-CNF formula
$\FF$, consider the weighted sum $F$
\[
    X = X(\gamma,\eta) = \sum_{\sigma} \gamma^{H(\sigma,{\FF})} \eta^{U(\sigma,{\FF})-u_0 m}
\enspace .
\]
Since $\gamma,\eta<1$ we see that in $X$ \chek{the} truth
assignments $\sigma$ for which $H(\sigma,F)>0$ or
$U(\sigma,{\FF})>u_0m$ are suppressed exponentially, whereas
\chek{the rest} are rewarded exponentially. Decreasing
$\gamma,\eta \in [0,1)$ makes this phenomenon more and more acute,
with the limiting case $\gamma, \eta = 0$ corresponding to a 0-1
weighting scheme (we adopt the convention $0^0\equiv 1$). Indeed,
applying the second moment method to $X$ with $\eta =0$
corresponds to the approach of~\cite{yuval} for the random $k$-SAT
threshold, where only satisfying assignments receive non-zero
weight $\gamma^{H(\sigma,F)}$. A key step in our analysis,
presented in Subsection \ref{tune}, is the tuning of the
parameters $\gamma,\eta$ to focus on truth assignments $\sigma$
for which (\ref{iopu}) and (\ref{colopu}) hold. Before doing that,
we establish the upper bound in Theorem~\ref{main}.

\subsection{The upper bound in Theorem \ref{main}}
This upper bound can be \chek{readily established by using the
entropic-form Chernoff bound for the Binomial} 
(see Lemma A.10 in \cite{AS} or Lemma 3.8 in \cite{DM}),
but it is more
informative to give a self-contained argument.
Recall the definition of $\rr_k(\cdot)$
from (\ref{defrr}).

\begin{lemma}\label{lemma:upper}
For all $k \geq 2$ and $\p \in (0,1]$, if $\q=1-p$ then
\begin{equation}\label{eq:upper}
r_k^{\up}(\p)\le \frac{2^k \log 2}
{\q \log \q -(2^{k}-\q)\log\Big(\frac{2^k-1}{2^k-q}\Big)}
\le {\rr}_k(\p) \enspace .
\end{equation}
\end{lemma}
\begin{proof}
The right hand inequality of (\ref{eq:upper}) follows from the
inequality $\log t \le t-1$ applied to $t=\frac{2^k-1}{2^k-q}$,
so we just need to verify the  left hand inequality. To do that,
write $u_0=2^k \q$. Let $\eta \in (0,1)$, and
observe that if  $\FF$ is $\p$-satisfiable, then $U(\sigma,{\FF})
\le u_0 m$ for some $\sigma$, whence
$$
X(1,\eta)=  \sum_{\sigma} \eta^{U(\sigma,{\FF})-u_0 m} \ge 1 \,.
$$
From (\ref{eq:mean}) in the next section we have that
\begin{equation}\label{eq:up2}
\PS[X(1,\eta) \ge 1] \le \ex[X(1,\eta)]=2^n \eta^{-\q rn 2^{-k}}
  \Big(1- (1-\eta)2^{-k}\Big)^{rn} \,.
\end{equation}
Thus, the probability of $\p$-satisfiability decays exponentially
in $n$ if the \chek{the $n$-th root of the RHS of} \eqref{eq:up2}
is strictly smaller than $1$.
Taking $\eta=\q (2^k-1)/(2^k-\q)$ 
yields the lemma.
\end{proof}

\subsection{Tuning parameters and truncation} \label{tune}

When $\eta>0$, attempting to apply the
second moment method to $X$ we encounter two major problems.

The first problem is that while $X>0$ implies satisfiability when
$\eta=0$, when $\eta >0$ having $X>0$ does not imply
$\p_0$-satisfiability:
in principle, $X$ could be positive due to the
contribution of assignments falsifying many more clauses than $u_0
m$. This necessitates restricting the sum defining $X$ to truth
assignments falsifying at most $u_0m+O(\sqrt{m})$ clauses, \ie
truncating $X$.

The second, more severe, problem is that with or without this
truncation, $\ex[X]^2/\ex[X^2]$ becomes exponentially small when
$r$ is only, roughly, half the (asymptotically optimal) lower
bound of Theorem~\ref{main}.  Rather counterintuitively, we will
be able to delay this explosion until $r$ is within $1-o(1)$ of
the upper bound by also removing from the sum those ``heroic"
truth assignments falsifying {\em fewer\/} than $u_0 m$ clauses.
This affords us much tighter control of pairs of assignments that
agree on nearly all variables, which turn out to be the dominant
contributors to $\ex[X^2]$ as we approach the upper bound. The
idea behind this sacrifice is motivated by Cramer's classical
``change of measure'' technique in large deviation theory. The
corresponding ``adaptive weighting" scheme requires an extremely
sharp asymptotic analysis, involving a number of rather miraculous
cancellations. Due to space limitations this analysis appears
entirely in the Appendix.

Specifically, for some fixed $A>0$ let
$${\cal S^*}= \{ \sigma \in \{0,1\}^n :
H(\sigma, {\FF}) \geq 0 \mbox{ \rm and } \, U(\sigma,\FF) \in [u_0
m ,u_0 m +A\sqrt{m}] \}  \, .
$$
Moreover, given $u_0$, let $\gamma_0,\eta_0$ be defined by
\begin{eqnarray}
1-\eta_0  & = & (1-\gamma_0^2)(1+\gamma_0^2)^{k-1} \nonumber\\
\label{eq:defzero} \\
u_0  & = & \frac{\eta_0}{(1+\gamma_0^2)^k-(1-\eta_0)} \nonumber
\enspace .
\end{eqnarray}
These two equations are designed so that the main contribution in
the sum defining $X$ comes from truth assignments  for which
(\ref{iopu}) and (\ref{colopu}) holds. The connection is made in
equations (\ref{EH}) and (\ref{EU}) in Section \ref{doobya}.

\chek{We define}
\[
X_*  = \sum_{\sigma \in {\cal S^*}}
\gamma_0^{H(\sigma,F)}\eta_0^{U(\sigma,F)-u_0 m} \enspace .
\]
Note that, by definition, when $X_*>0$ at least one truth
assignment must falsify at most $u_0 m + A\sqrt{m}$ clauses. Thus,
if for a given $\p_0$ we can prove that there exists a constant
$D>0$ such that $\ex[X_*^2] < D \times \ex[X_*]^2$ then, by
Corollary~\ref{booster}, it follows that $F_k(n,rn)$ is \whp\
$\p$-satisfiable for all $\p<\p_0$.

Bounding the second moment of $X_*$ will be accomplished in the
following lemmata. For $\alpha \in [0,1]$, let
\begin{flalign}\label{eq:deff}
\mbox{} \!\!\!\!\! f(\alpha,\gamma,\eta) =\eta^{-2u_0}\left[
    \left(\a \left(\frac{\gamma^2+\gamma^{-2}}{2}\right)+1-\a
    \right)^k
    - 2 (1-\eta)\left(\frac{\a\gamma^{-2}+(1-\a)}{2}\right)^k
    + (1-\eta)^2 \left(\frac{\a\gamma^{-2}}{2}\right)^k\right]
\end{flalign}
and
\begin{equation}
 g_r(\alpha,\gamma,\eta) = \frac{f(\alpha,\gamma,\eta)^r}{\alpha^{\alpha}(1-\alpha)^{1-\alpha}}\enspace .
\end{equation}\medskip
In all of the following lemmata $k \geq 2$ is a fixed integer and
$r>0$.
\begin{itemize}
\item
Lemma~\ref{middle} with $\gamma=\gamma_0$ and $\eta=\eta_0$ gives
us $\ex[X(\gamma_0,\eta_0)]^2$ which is $\ex[X_*]^2$ but for the
truncation.
\item
Lemma~\ref{doobya2} asserts that for every value of $u_0$,
$\ex[X_*]$ is a constant fraction of $\ex[X(\gamma_0,\eta_0)]$.
Thus, combined with Lemma~\ref{middle}, it gives us $\ex[X_*]^2$
up to a constant factor (which is all we need).
\item
Lemma~\ref{ola} expresses $\ex[X_*^2]$ as a sum with $n+1$ terms,
the $z$-th term capturing the contribution of the $2^n
\binom{n}{z}$ pairs of truth assignments with overlap $z$. The
contribution of each such pair is then bounded by
$f(z/n,\gamma,\eta)^{rn}$ where $\gamma,\eta$ are allowed to {\em
depend on $z$}, subject only to $\gamma \geq \gamma_0$ and $\eta
\geq \eta_0$ respectively. In other words, Lemma~\ref{ola} allows
us to adapt $\gamma$ and $\eta$ to $\alpha$, which is crucial when
$\p <1$.
\item
Lemma~\ref{ksekola} is based on the fact that for any ``smooth"
choice of sequences $\gamma(z),\eta(z)$, the sum in
Lemma~\ref{ola} will be dominated by the contribution of the
$\Theta(n^{1/2})$ terms around the maximum term. Specifically, if
$\chi,\omega$ express our adaptive scheme for $\gamma,\eta$, then
we can use the Laplace method to get that the maximum of
$g_r(\alpha,\chi(\alpha),\omega(\alpha))$ over $\alpha \in (0,1)$,
characterizes the sum in Lemma~\ref{ola} up to a constant factor.
\end{itemize}

\begin{lemma}\label{middle}
For every $u_0,\gamma,\eta \in [0,1)$,
\[
\ex[X]^2 = \Bigl(2\,g_r(1/2,\gamma,\eta)\Bigr)^n \enspace .
\]
\end{lemma}

\begin{lemma}\label{doobya2}
For every $u_0$, there exists $\theta = \theta(k,A)>0$ such that
as $n \to \infty$,
\[
\frac{\ex[X_*]}{\ex [X(\gamma_0,\eta_0)]} \to \theta \enspace .
\]
\end{lemma}

\begin{lemma}\label{ola}
Let $\gamma(z),\eta(z)$ be arbitrary sequences such that
$\gamma(z) \geq \gamma_0$ and $\eta(z) \geq \eta(0)$ for every $0
\leq z \leq n$. Then, for every $u_0$,
\[
\ex[X_*^2] \leq 2^n \sum_{z=0}^n \binom{n}{z}
f(z/n,\gamma(z),\eta(z))^{rn} \enspace .
\]
\end{lemma}

\begin{lemma}\label{ksekola}
Let $\chi : [0,1] \rightarrow [\gamma_0,1)$ and $\omega : [0,1]
\rightarrow [\eta_0,1)$ be arbitrary piecewise-smooth functions
and let $g_r(\alpha) = g_r(\alpha,\chi(\alpha),\omega(\alpha))$.
If there exists $\amax \in (0,1)$ such that $g_r(\amax) \equiv
\gmax > g_r(\alpha)$ for all $\alpha \neq \amax$, and
$g''_r(\amax)< 0$, then there exists a constant $D =
D_{\chi,\omega}(k,r,u_0)>0$ such that for all sufficiently large
$n$
\[
\ex[X_*^2] < D \times {
    \Bigl(2 \,
            \gmax
    \Bigr)^n}\enspace .
\]
\end{lemma}

Combining Lemmata~\ref{middle}--\ref{ksekola} we  see that if for
a given $u_0$ and $r$ there exist $\chi,\omega$ such that for all
$\alpha \neq 1/2$
\begin{equation}\label{masterp}
g_{r}\Bigl(1/2,\gamma_0,\eta_0\Bigr) >
g_{r}\Bigl(\alpha,\chi(\alpha),\omega(\alpha)\Bigr)
\end{equation}
then $\ex[X_*^2] < D \theta^{-2} \times \ex[X_*]^2$, yielding the
desired conclusion $\ex[X_*^2] = O(\ex[X_*]^2)$.
\bigskip

Indeed, to prove  Theorem~\ref{main} we will show that for every
$\p \in (0,1]$ and for the stated $r = r(\p)$, there exist
functions $\chi,\omega$ for which~\eqref{masterp} holds. To
simplify the asymptotic analysis, we use the crudest possible such
functions, paying the price of this simplicity in the value of
$k_0$ in Proposition~\ref{prop:main} below. We note that by
choosing a more refined (and more cumbersome) adaptation of
$\gamma,\eta$ to $\alpha$ this value can be improved greatly.
Moreover, we emphasize that for any fixed value of $k$, one can
get a sharper lower bound (such as those reported in the
Introduction) by partitioning $[0,1]$ to a large number of
intervals and numerically finding a good value of $\gamma,\eta$
for each one. We discuss this point further in
Section~\ref{finitek}. Finally, we note that general large
deviations considerations imply that for every $k$ and $\p$, the
condition~\eqref{masterp} is sharp for our method. That is, no
better lower bound can be derived by considering balanced
assignments and, in fact, by any argument that classifies
assignments according to their number of satisfied literal
occurrences in the formula.
\begin{definition}\label{skolio}
Let $\q = 1-\p_0 = u_0 2^k$ and let
\begin{eqnarray}\label{eq:defrk}
\ttk=\frac{2^k\ln 2}{1-\q+\q\log
\q}\left(1-20k2^{-k\varphi(\q)}\right)\quad \mathrm{where}\quad
\varphi(\q)=\frac{(1-\sqrt{\q})^2}{1-\q+\q\log \q} \enspace .
\end{eqnarray}
\end{definition}
Theorem~\ref{main} will follow from the following
Proposition.
\begin{proposition}\label{prop:main}
Let
\begin{eqnarray}\label{eq:defE}
G_r(\alpha)=\left\{\begin{array}{ll} g_r(\alpha,\gamma_0,\eta_0)
&\mbox{ if }
\alpha \in \left[\frac{3\log k}{k}, 1-\frac{3\log k}{k}\right]\\
\\
g_r(\alpha, \sqrt{\gamma_0},\sqrt{\eta_0}) & \mbox{ otherwise.}
\end{array}\right.
\end{eqnarray}
For all $k\ge k_0$, if\/ $r\le \ttk$ then $G_r''(1/2)<0$ and
$G_r(1/2)>G_r(\alpha)$ for all $\alpha\neq 1/2$.
\end{proposition}

\bigskip

The proof of Proposition~\ref{prop:main}, itself, will be
decomposed into three lemmata of increasing difficulty. The first
lemma holds for any $\gamma,\eta$ and reduces the proof to the
case $\alpha \geq 1/2$. The second lemma reflects the behavior of
$f$ (and thus $g_r$) around $\a = 1/2$, motivating the judicious
choice $\eta=\eta_0$ and $\gamma=\gamma_0$ for $G_r$. The third
lemma deals with $\alpha$ near 1. That case needs a lot more work
in order to handle the unique local maximum of $g_r$ in that
region. The condition $r \leq \ttk$ and the change to
$\gamma=\sqrt{\gamma_0}, \eta=\sqrt{\eta_0}$ aims precisely at
keeping the value of $g_r$ at this other local maximum smaller
than $g_r(1/2,\gamma_0,\eta_0)$.

\begin{lemma}
\label{lem:biggerthanhalf} For every $0< x\le \frac12$,
$G_r(1/2+x)>G_r(1/2-x)$.
\end{lemma}

\begin{lemma}
\label{lem:decreasing} For all $k\ge k_0$, if\/ $r\le
\frac{2^k\log 2}{1-\q+\q\log \q}$ then $G_r''(1/2)<0$ and $G_r$ is
strictly decreasing on $\left[\frac{1}{2},1-\frac{3\log
k}{k}\right]$.
\end{lemma}

\begin{lemma}
\label{lem:bigalpha} For all $k \geq k_0$, if  $r\le \ttk$ then
for every $\alpha\in \left[1-\frac{3\log k}{k},1\right]$,
$G_r(1/2)>G_r(\alpha)$.
\end{lemma}

In the following sections we prove
Lemmata~\ref{middle}--\ref{ksekola}, while
Lemmata~\ref{lem:biggerthanhalf}--\ref{lem:bigalpha} are proven in
the appendix. Before delving into the probabilistic calculations
involved in proving Lemmata~\ref{middle}--\ref{ola} a couple of
remarks are in order.\medskip

\noindent {\bf Relationship to other $k$-CNF models:} Recall that
the $m$ clauses of $F_k(n,m)$ are chosen independently with
replacement among the $(2n)^k$ possibilities. Thus, the $m$
clauses $\{c_i\}_{i=1}^m$ are i.i.d.\ random variables, each $c_i$
being the conjunction of $k$ i.i.d.\ random variables
$\{\ell_{ij}\}_{j=1}^k$, each $\ell_{ij}$ being a uniformly random
literal. This viewpoint of the formula as a sequence of $km$
i.i.d.\ random literals will be very handy for our calculations.

Clearly, in this model some clauses might be improper, \ie they
might contain repeated and/or contradictory literals. At the same
time, though, observe that the probability that any given clause
is improper is smaller than $k^2/n$ and, moreover, the proper
clauses are uniformly selected among all such clauses. Therefore
\whp\ the number of improper clauses is $o(n)$ implying that if
for a given $r$, $F_k(n,rn)$ is $\p$-satisfiable \whp\ then for
$m=rn-o(n)$, the same is true in the model where we only select
among proper clauses. The issue of selecting clauses without
replacement is completely analogous as \whp\ there are $o(n)$
clauses that contain the same $k$ variables as some other clause.
\medskip

\noindent {\bf Notation:} In the ensuing probabilistic
calculations it will be convenient to write $\sigma \viol F$ to
denote that the truth assignment $\sigma$ violates the formula $F$
where $F$ can be a literal, a clause, or an entire CNF.

\section{The first moment and proof of Lemma~\ref{middle}}

By linearity of expectation and since the $m = rn$ clauses
$c_1,c_2,\ldots,c_m$ are chosen independently we have
\begin{eqnarray}
\eta^{u_0 m}  \, \ex[X]   & = &  \ex \left[\sum_{\sigma}  \gamma^{H(\sigma,{\FF})}  \eta^{U(\sigma,{\FF})}\right] \nonumber \\
            & = &  \sum_{\sigma}   \ex\left[\prod_{c_i}
                            \gamma^{H(\sigma,c_i)}  \eta^{U(\sigma,c_i)}
                    \right] \nonumber \\
            & = & \sum_{\sigma}   \prod_{c_i}\ex
                    \left[
                            \gamma^{H(\sigma,c_i)}  \eta^{U(\sigma,c_i)}
                    \right] \enspace . \label{eq:sumonday}
\end{eqnarray}
%
%

Observe now that since the clauses are identically distributed, by
symmetry, it suffices to consider the expectation
in~\eqref{eq:sumonday} for a single random clause $c= \ell_1 \vee
\cdots \vee \ell_k$ and a fixed truth assignment $\sigma$.
Moreover, observe that if we write $\gamma^H \eta^U$ as $\gamma^H
+ \gamma^H(\eta^U-1)$ we see that the second expression is
non-zero only when $U>0$, \ie when $c$ is violated by $\sigma$.
So, since the literals $\ell_1,\ldots,\ell_k$ are i.i.d.\ we get
\begin{eqnarray}
    \ex\left[\gamma^{H(\sigma,c)} \eta^{U(\sigma,c)}\right] & = &
    \ex\left[\gamma^{H(\sigma,c)} - \gamma^{H(\sigma,c)}
\left(1-\eta^{U(\sigma,c)}\right)\right]\nonumber  \\
    & = & \ex\left[ \gamma^{H(\sigma,c)}\right] -
 \ex\left[\gamma^{-k} (1-\eta) \, \one_{\sigma \viol c}\right] \nonumber \\
& = &     \ex\left[\prod_{\ell_i}\gamma^{H(\sigma,\ell_i)}\right]
- 2^{-k} \gamma^{-k} (1-\eta) \nonumber \\
& = &  \prod_{\ell_i}\ex\left[\gamma^{H(\sigma,\ell_i)}\right]
- 2^{-k} \gamma^{-k} (1-\eta) \nonumber \\
& = & \left(\frac{\gamma+\gamma^{-1}}{2}\right)^k - (2\gamma)^{-k}
(1-\eta)
\nonumber \\
& \equiv & \Z(\gamma,\eta) \label{eq:defZ} \enspace .
\end{eqnarray}
Thus,
\begin{equation} \label{eq:mean}
\ex[X] = \eta^{-u_0 rn} \, 2^n \, \Z(\gamma,\eta)^{rn} \,.
\end{equation}
Observe now that
\[
\left(\eta^{-u_0}\Z(\gamma,\eta)\right)^2 = f(1/2,\gamma,\eta)
\enspace .
\]
Therefore,
\[
\ex[X]^2 = \left(\eta^{-u_0 rn} \, 2^n \,
\Z(\gamma,\eta)^{rn}\right)^2 =
 \left[\left(\eta^{-u_0 r} \, 2 \,
\Z(\gamma,\eta)^r\right)^2\right]^n = \left[4 \,
f(1/2,\gamma,\eta)^r\right]^{n} = \left[2 \,
g_r(1/2,\gamma,\eta)\right]^n \enspace .
\]

\section{Proof of Lemma~\ref{doobya2}} \label{doobya}

By linearity of expectation, it suffices to prove that there
exists some $\theta=\theta(k,A)>0$ such that for the values of
$\gamma_0, \eta_0$ satisfying \eqref{eq:defzero} and every truth
assignment $\sigma$, we have
\begin{equation} \label{favor}
    \frac{\ex\left[\gamma_0^{H(\sigma,\FF)}\eta_0^{U(\sigma,\FF)}
\one_{\sigma \in {\cal S^*({\FF})}}\right ]}
         {\ex \left[\gamma_0^{H(\sigma,\FF)}\eta_0^{U(\sigma,\FF)}\right]}
 \rightarrow \theta \enspace .
\end{equation}

Recalling that formulas in our model are sequences of i.i.d.\
random literals $\ell_1,\ldots,\ell_{km}$, let $\PS(\cdot)$ denote
the probability assigned by our distribution to any such sequence,
\ie $(2n)^{-km}$. Now, fix any truth assignment $\sigma$ and
consider an auxiliary distribution $\Pa$ on $k$-CNF formulas where
the $m$ clauses $c_1,\ldots,c_m$ are again i.i.d.\ among all
$(2n)^k$ clauses, but where now for any fixed clause $\fc$
\begin{equation} \label{doob2}
\Pa(\rc_i=\fc) =
\frac{\gamma_0^{H(\sigma,\fc)}\eta_0^{U(\sigma,\fc)}\PS(\fc)}
{\Z(\gamma_0,\eta_0)} \enspace ,
\end{equation}
where
\begin{equation}\label{defZ2}
\Z(\gamma_0,\eta_0) =
\ex\left[\gamma_0^{H(\sigma,\rc)}\eta_0^{U(\sigma,\rc)} \right]
\enspace ,
\end{equation}
was defined in~(\ref{eq:defZ}). (Since each fixed clause $\fc$
receives probability proportional to
$\gamma_0^{H(\sigma,\fc)}\eta_0^{U(\sigma,\fc)}$, indeed
$Z(\gamma_0,\eta_0)$ provides the correct normalization to a
probability distribution.) So, whereas under $\PS(\cdot)$ every
$k$-CNF formula ${\FF}$ with $m$ clauses had the same probability
$\PS(F)=(2n)^{-km}$, under $\Pa$ its probability is
\begin{equation} \label{doob3} \Pa({\FF}) =
\frac{\gamma_0^{H(\sigma,{\FF})}\eta_0^{U(\sigma,{\FF})}
\PS({\FF})} {\Z(\gamma_0,\eta_0)^m} \enspace .
\end{equation}

Let $\Ea$ be the expectation operator corresponding to $\Pa$. A
calculation similar to that leading to (\ref{eq:defZ}), adding the
equal contributions from the $k$ literals, gives that for a single
random clause $c$
\begin{equation} \label{EH}
\Z(\gamma_0,\eta_0)\Ea[H(\sigma,\rc)]= k(\gamma_0-\gamma_0^{-1})
\left(\frac{\gamma_0+\gamma_0^{-1}}{2}\right)^{k-1}
+k(2\gamma_0)^{-k} (1-\eta_0) \enspace .
\end{equation}
Moreover,
\begin{equation} \label{EU}
\Z(\gamma_0,\eta_0)\Ea[U(\sigma,\rc)]= (2\gamma_0)^k \eta_0 \,.
\end{equation}
Thus~\eqref{eq:defzero} ensures that $\Ea[H(\sigma,\rc)] = 0$ and
also that $\Ea[U(\sigma,\rc)-u_0] = 0$.\medskip

Next, we apply the multivariate central limit theorem (see, e.g.\
\cite{pollard}, page 182) to the i.i.d.\ mean-zero random vectors
$\Big(H(\sigma,\rc_i),U(\sigma,\rc_i)-u_0\Big)$ for
$i=1,\ldots,m$. Observe that, since $k \ge 2$,  the common law of
these random vectors is not supported on a line. We deduce that as
$n \rightarrow \infty$
$$
\Pa[\sigma \in {\cal S}^*(\FF)]=\Pa\Big[H(\sigma,{\FF}) \geq 0 \mbox{ \rm
and } U(\sigma,\FF) \in [m u_0 ,m u_0+A\sqrt{m}]\Big]
 \rightarrow \theta(k,A)>0 \enspace .
$$
Here, the right hand side is the probability that a certain
nondegenerate bivariate normal law assigns to a certain open set.
Its exact value is unimportant for our purpose. By (\ref{doob3}),
this is equivalent to (\ref{favor}).

\section{Proof of Lemma~\ref{ola}}

Linearity of expectation implies
\begin{eqnarray}
\eta_0^{2u_0 m} \, \ex [X_*^2]     & = & \ex \left[\left(\sum_{\sigma}  \gamma_0^{H(\sigma,{\FF})}  \eta_0^{U(\sigma,{\FF})} \, \one_{\sigma \in {\cal S^*}({\FF})}\right)^2 \right] \nonumber \\
            & = &   \sum_{\sigma,\tau} \,\ex \left[\gamma_0^{H(\sigma,{\FF})+H(\tau,{\FF})}
            \eta_0^{U(\sigma,{\FF})+U(\tau,{\FF})}\, \one_{\sigma,\tau \in {\cal S^*}({\FF})}
        \right] \label{eq:doldrum} \enspace .
\end{eqnarray}
Observe now that since $\sigma \in {\cal S}^*$ implies
$H(\sigma,F) \geq 0$ and $U(\sigma,F) \geq u_0 m$, we get that for
every pair $\sigma,\tau$ and any $\gamma \ge \gamma_0$ and $\eta
\ge \eta_0$,
\begin{eqnarray}
    \ex \left[\gamma_0^{H(\sigma,{\FF})+H(\tau,{\FF})}
            \eta_0^{U(\sigma,{\FF})+U(\tau,{\FF})}\, \one_{\sigma,\tau \in {\cal S^*}({\FF})}
        \right]    &\le&
\ex\left[\gamma^{H(\sigma,{\FF})+H(\tau,{\FF})}
\eta^{U(\sigma,{\FF})+U(\tau,{\FF})} \, \one_{\sigma,\tau \in
{\cal S^*}({\FF})}\right] \nonumber \\
   &\le&
 \ex\left[\gamma^{H(\sigma,{\FF})+H(\tau,{\FF})}
\eta^{U(\sigma,{\FF})+U(\tau,{\FF})}\right] \enspace .
\label{eq:choice}
\end{eqnarray}
In other words, when using the right hand side
of~\eqref{eq:choice} to bound each term of the sum
in~\eqref{eq:doldrum}, we are allowed to {\em adapt\/} the value
of $\gamma$ and $\eta$ to the pair $\sigma,\tau$, the only
restrictions being $\gamma \geq \gamma_0$ and $\eta \geq \eta_0$.
This is a crucial point and we will exploit it heavily when
bounding the contribution of pairs with large overlap.
\medskip

To estimate the right hand side of~\eqref{eq:choice} for any pair
$\sigma,\tau$ we first observe that since the $m$ clauses
$c_1,c_2,\ldots,c_m$ are i.i.d., letting $c$ be a single random
clause we have
\begin{eqnarray}
\ex\left[\gamma^{H(\sigma,{\FF})+H(\tau,{\FF})}
\eta^{U(\sigma,{\FF})+U(\tau,{\FF})}\right]        & = &
        \ex\left[\prod_{c_i} {\gamma^{H(\sigma,c_i)+H(\tau,c_i)}
\eta^{U(\sigma,c_i)+U(\tau,c_i)}}
                    \right] \nonumber \\
            & = &  \prod_{c_i}\ex
                    \left[{\gamma^{H(\sigma,c_i)+H(\tau,c_i)}
\eta^{U(\sigma,c_i)+U(\tau,c_i)}}
                    \right]\nonumber \\
                    & = &
                    \biggl(\ex \left[{\gamma^{H(\sigma,c)+H(\tau,c)}
\eta^{U(\sigma,c)+U(\tau,c)}}
                    \right]\biggr)^m \label{eq:monday} \enspace .
\end{eqnarray}

Next, we observe that for every pair $\sigma,\tau$, by symmetry,
the expectation in~\eqref{eq:monday} depends only on the number of
variables to which $\sigma,\tau$ assign the same value. So, let
$\sigma,\tau$ be any pair of truth assignments that agree on
exactly $z = \alpha n$ variables, \ie have overlap $z$. By first
rewriting (again) $\gamma^H \eta^U$ as $\gamma^H +
\gamma^H(\eta^U-1)$ and then observing that $\eta^{U(\tau,c)}$ is
distributed identically with $\eta^{U(\sigma,c)}$ we  get
\begin{align}
\lefteqn{\ex
    \biggl[
              \gamma^{H(\sigma,{c})}\eta^{H(\sigma,{c})} \gamma^{H(\tau,{c})}  \eta^{H(\tau,{c})}
    \biggr]} \nonumber \\
\mbox{} \quad \mbox{} & = \ex
    \Biggl[
              \biggr(\gamma^{H(\sigma,c)}-\gamma^{H(\sigma,c)}\left(1-\eta^{U(\sigma,c)}\right)\biggr)
              \biggr(\gamma^{H(\tau,c)}-\gamma^{H(\tau,c)}\left(1-\eta^{U(\tau,c)}\right)\biggr)
    \Biggr]
\nonumber \\
& =        \ex\biggl[\gamma^{H(\sigma,c)+H(\tau,c)}\biggr]
            - 2 \,
            \ex\biggl[\gamma^{H(\sigma,c)+H(\tau,c)}\left(1-\eta^{U(\sigma,c)}\right)\biggr]+
\ex\biggl[\gamma^{H(\sigma,c)+H(\tau,c)}\left(1-\eta^{U(\sigma,c)}\right)\left(1-\eta^{U(\tau,c)}\right)\biggr]
\nonumber
            \\
& =        \ex\left[\gamma^{H(\sigma,c)+H(\tau,c)}\right]
            - 2(1-\eta)
            \ex\left[\gamma^{H(\sigma,c)+H(\tau,c)} \one_{\sigma \viol c}\right]+
            2^{-k}\alpha^k \gamma^{-2k} (1-\eta)^2   \enspace .
            \label{eq:merge}
\end{align}

Now, to estimate~\eqref{eq:merge} we note that since the literals
$\ell_1,\ell_2,\ldots \ell_k$ comprising $c$ are i.i.d.\ we have
\begin{align*}
    \ex\left[\gamma^{H(\sigma,c)+H(\tau,c)}
    \right] & =
    \ex\left[
            \prod_i \gamma^{H(\sigma,\ell_i)+H(\tau,\ell_i)}
        \right]
      =
        \prod_i \ex\left[\gamma^{H(\sigma,\ell_i)+H(\tau,\ell_i)}
        \right]
        = \left(\a \left(\frac{\gamma^2+\gamma^{-2}}{2}\right)+1-\a
        \right)^k \nonumber \\
\intertext{and, similarly,}
\ex\left[\gamma^{H(\sigma,c)+H(\tau,c)} \one_{\sigma \viol c}
    \right] & =
    \ex\left[
            \prod_i \gamma^{H(\sigma,\ell_i)+H(\tau,\ell_i)} \one_{\sigma
            \viol \ell_i}
        \right]
     =
        \prod_i \ex\left[\gamma^{H(\sigma,\ell_i)+H(\tau,\ell_i)} \one_{\sigma \viol \ell_i}
        \right]  = \left(\frac{\a\gamma^{-2}+(1-\a)}{2}
        \right)^k .
\end{align*}
Substituting these last two equations in~\eqref{eq:merge} we get
\begin{align} \lefteqn{\eta ^{-2u_0}\,\ex
    \left[
              \gamma^{H(\sigma,{c})}\eta^{H(\sigma,{c})} \gamma^{H(\tau,{c})}  \eta^{H(\tau,{c})}
    \right]} \nonumber \\
& = \nonumber \eta^{-2u_0}\left[
    \left(\a \left(\frac{\gamma^2+\gamma^{-2}}{2}\right)+1-\a
    \right)^k
    - 2 (1-\eta)\left(\frac{\a\gamma^{-2}+(1-\a)}{2}\right)^k
    + (1-\eta)^2 \left(\frac{\a\gamma^{-2}}{2}\right)^k\right] \\
& =  f(\alpha,\gamma,\eta) \enspace . \label{eq:gotcha}
\end{align}

So, in conclusion, since the number of ordered pairs with overlap
$z$ is $ 2^n \,\binom{n}{z}$ we get that
\begin{equation}
 \ex[X_*^2] \leq 2^n \sum_{z=0}^n  \binom{n}{z} f(z/n,\gamma(z),\eta(z))^{m}
 \enspace ,
\label{eq:naesum}
\end{equation}
for any set of choices for $\gamma(z),\eta(z)$ such that
$\gamma(z) \geq \gamma_0$ and $\eta(z)\geq \eta_0$ for all $0\leq
z \leq n$.\bigskip

\subsection{Proof of Lemma~\ref{ksekola}}

If $\chi : [0,1] \rightarrow [\gamma_0,1]$ and $\omega : [0,1]
\rightarrow [\gamma_0,1]$ are piecewise smooth, then from the
definition of $f$ we see that
$f(\alpha,\chi(\alpha),\omega(\alpha))$ is also piecewise smooth.
Thus, we can decompose the sum in~\eqref{eq:naesum} into a fixed
number of sums such that $f(\alpha,\chi(\alpha),\omega(\alpha))$
is smooth in the range of each sum. To bound each such sum, then,
we use the following lemma whose proof is implied by the proof of
Lemma~2 in~\cite{naesat} (that lemma is stated with the
requirement that $f$ is analytic, a condition not needed for the
proof; in fact, it suffices for $f$ to only be twice
differentiable.) The idea is that each of these sums is dominated
by the contribution of $\Theta(n^{1/2})$ terms around the maximum
term. Since the number of sums is finite the lemma follows.
\begin{lemma} Let $\phi$ be any
real, positive, twice-differentiable function on $[0,1]$ and let
$$
S_n = \sum_{z=0}^n \binom{n}{z} \,\phi(z/n)^n \enspace .
$$
Letting $0^0 \equiv 1$, define $g$ on $[0,1]$ as
\[
g(\alpha) = \frac{\phi(\alpha)}
                 {\alpha^\alpha \,(1-\alpha)^{1-\alpha}} \enspace
                 .
\]
If there exists $\amax \in (0,1)$ such that $g(\amax) \equiv \gmax
> g(\alpha)$ for all $\alpha \neq \amax$, and $g''(\amax) < 0$,
then there exist constants $B, C > 0$ such that for all
sufficiently large $n$
\[
B \times \gmax^n \,\le\, S_n \,\le\, C \times \gmax^n \enspace .
\]
\end{lemma}

\section{Bounds for finite {\large{$k$}}}\label{finitek}

As mentioned in Section~\ref{outline}, for small values of $k$ the
simple adaptation scheme of Proposition~\ref{prop:main} does not
yield the best possible lower bound for $\p$-satisfiability
afforded by our method. For that, one has to use a significantly
more refined adaptation of $\gamma,\eta$ with respect to $\alpha$.
Our lower bounds reported in Figure 1 are, indeed, the result of
performing such optimization of $\gamma,\eta$ numerically (for
both the upper bound plots and the plots of the lower bound
from~\cite{sork} we used the explicit formulas).

Specifically, to create the plots of the lower bounds we computed
a lower bound for 100 equally spaced values of $p$ on the
horizontal axis (and then had Maple's~\cite{maple} plotting
function ``connect the dots"). For each of these values of $p$, to
prove the corresponding lower bound for $r$ we had to establish
that there exist a choice of functions $\chi,\omega$ as in
Lemma~\ref{ksekola} such that for all $\alpha \in (1/2,1]$ we have
$g_{r}(1/2,\gamma_0,\eta_0) >
g_{r}(\alpha,\chi(\alpha),\omega(\alpha))$. To that end, we
partitioned $(1/2,1]$ to 10,000 points and for each such point we
searched for values of $\gamma \geq \gamma_0$ and $\eta \geq
\eta_0$ such that this condition holds with a bit of room. (For
$k>4$ we solved~\eqref{eq:defzero}, defining $\gamma_0$ and
$\eta_0$, numerically to 10 digits of accuracy. For the
optimization we exploited convexity to speed up the search.)
Having determined such values, we (implicitly) extended the
functions $\chi,\omega$ to all $(1/2,1]$ by assigning to every
not-chosen point the value at the nearest chosen point. Finally,
we computed a (crude) upper bound on the derivative of $g_r$ with
respect to $\alpha$ in $(1/2,1]$. This bound on the derivative,
along with our room factor, then implied that for every point that
we did not check, the value of $g_r$ was sufficiently close to its
value at the corresponding chosen point to also be dominated by
$g_{r}(1/2,\gamma_0,\eta_0)$.

\section*{Acknowledgements}
We thank Cris Moore for helpful conversations in the early stages
of this work.

\bibliographystyle{amsalpha}
\bibliography{stoc,extra,theory}

\newpage

\appendix

\section{Building up an arsenal}

In this section we collect some basic inequalities and identities
that we will use in the proofs of Lemmas \ref{lem:biggerthanhalf},
\ref{lem:decreasing} and \ref{lem:bigalpha}. For readability, in
this Appendix we have replaced $q$ of Definition~\ref{skolio} with
the letter $y$.

Plugging in the definition of $\eta_0$ from (\ref{eq:defzero})
into the definition of $f$ we get
\begin{align*} \lefteqn{f(\alpha,\gamma_0,\eta_0)} \nonumber \\
&
=\eta_0^{-{y}/{2^{k-1}}}\left\{\left(1-\alpha+\alpha\frac{\gamma_0^2+\gamma_0^{-2}}{2}\right)^k-
\frac{(1-\gamma_0^2)(1+\gamma_0^2)^{k-1}}{2^{k-1}}
\left(\alpha\gamma_0^{-2}+1-\alpha\right)^k+
\frac{\alpha^k(1-\gamma_0^2)^2(1+\gamma_0^2)^{2k-2}}{2^k\gamma_0^{2k}}\right\}.
\end{align*}
For some parts of the ensuing calculations, it will be to
convenient to use the following normalizations of
$f(\alpha,\gamma_0,\eta_0)$ and $g_r(\alpha,\gamma_0,\eta_0)$
denoted as $f_0$ and $g_0$ respectively
\begin{eqnarray}\label{eq:normalization}
f_0(\alpha)=2^{2k}\gamma_0^{2k}\eta_0^{y/2^{k-1}}f(\alpha,\gamma_0,\eta_0)\quad
\mathrm{and}\quad
g_0(\alpha)=2^{2kr}\gamma_0^{2kr}\eta_0^{yr/2^{k-1}}g_r(\alpha,\gamma_0,\eta_0).
\end{eqnarray}
We will also write $\e_0=1-\gamma_0^2$. With this notation, we
have the following formula for $f_0$, which holds for every
$x\in[-1/2,1/2]$
\begin{eqnarray}\label{eq:f0}
f_0\left(\frac{1}{2}+x\right)=[2x{\e_0}^2+(2-{\e_0})^2]^k-2{\e_0}(2-{\e_0})^{k-1}[2-{\e_0}+2x{\e_0}]^k+{\e_0}^2(2-{\e_0})^{2k-2}(1+2x)^k
\enspace ,
\end{eqnarray}
In particular,
\begin{eqnarray}\label{eq:f1/2}
f_0\left(\frac{1}{2}\right)=(2-{\e_0})^{2k}-2{\e_0}(2-{\e_0})^{2k-1}+{\e_0}^2(2-{\e_0})^{2k-2}
=4(1-{\e_0})^2(2-{\e_0})^{2k-2} \enspace .
\end{eqnarray}

The function $y\mapsto 1-y+y\log y$, defined on $[0,1]$, appears
throughout our analysis. The following inequalities, valid for all
$y\in [0,1]$, will be used
\begin{eqnarray}\label{eq:stupid}
\frac{(1-y)^2}{2}\le 1-y+y\log y\le (1-y)^2.
\end{eqnarray}
The right-hand inequality follows from the estimate $\log y\le
y-1$. The left-hand inequality follows from integrating this
estimate as follows
$$
-1+y-y\log y=\int_y^1\log x dx\le
\int_y^1(1-x)dx=-\frac{(1-y)^2}{2}.
$$

\medskip
We end this section by providing some estimates for the values of
${\e_0}$ and $\eta_0$
\begin{fact}\label{fact:etae} For all sufficiently large $k$ ,
\begin{eqnarray}\label{eq:e}
\frac{2(1-y)}{2^{k}-k-1}-\frac{4k(1-y)^2}{2^{2k}}\;\le\;
{\e_0}\;\le \;\frac{2(1-y)}{2^{k}-k-1} \enspace ,
\end{eqnarray}
and
\begin{eqnarray}\label{eq:eta}
\eta_0 \;\le \;
\min\left\{y,y-\frac{(k+1)(1-y)}{2^k-k-1}+\frac{4k(1-y)^2}{2^k}\right\}
\enspace .
\end{eqnarray}
\end{fact}

\subsection{Proof of Fact \ref{fact:etae}}

By the first equation in (\ref{eq:defzero}) we have that
$\eta_0=1-{\e_0}(2-{\e_0})^{k-1}$. Plugging this into the second
equation in (\ref{eq:defzero}) we find that
\begin{eqnarray}\label{eq:identityy}
\frac{y}{2^k}=\frac{1-{\e_0}(2-{\e_0})^{k-1}}{(2-{\e_0})^k-{\e_0}(2-{\e_0})^{k-1}}
= \frac{1-{\e_0}(2-{\e_0})^{k-1}}{2(1-{\e_0})(2-{\e_0})^{k-1}}
=\frac{1}{2}\left[1-\sum_{j=1}^{k-1}\frac{1}{(2-{\e_0})^{j}}\right].
\end{eqnarray}
Hence, if we denote
\begin{eqnarray}\label{eq:defpsi}
\psi(t)=\sum_{j=1}^{k-1}\frac{1}{(2-t)^{j}}=\frac{(2-t)^{k-1}-1}{(1-t)(2-t)^{k-1}}
\enspace ,
\end{eqnarray}
then we require that $\psi({\e_0})=1-\frac{y}{2^{k-1}}$. We
collect below some useful properties of $\psi$.

\begin{lemma}\label{lem:psi} $\psi$ is increasing on $[0,1]$.
Furthermore, for $k$ large enough and every $t\le t\le 1/(2k)$,
\begin{eqnarray}\label{eq:ineqpsi}
1-\frac{1}{2^{k-1}}+t-\frac{(k+1)t}{2^k}+\frac{t^2}{2}\le
\psi(t)\le  1-\frac{1}{2^{k-1}}+t-\frac{(k+1)t}{2^k}+2t^2
\enspace.
\end{eqnarray}

\end{lemma}

\begin{proof} The fact that $\psi$ is increasing follows
immediately from the first formula in (\ref{eq:defpsi}). To prove
the inequalities in (\ref{eq:ineqpsi}), observe that
$$
\psi(t)=\frac{\left(1-\frac{t}{2}\right)^{k-1}-\frac{1}{2^{k-1}}}{(1-t)\left(1-\frac{t}{2}\right)^{k-1}}=
\frac{1}{1-t}-\frac{1}{2^{k-1}(1-t)\left(1-\frac{t}{2}\right)^{k-1}}
\enspace .
$$
To estimate $\psi$ from below, we use the inequalities $1+a+a^2\le
1/(1-a)\le 1+a+2a^2$ and $(1-a)^{k-1}\ge1-(k-1)a$, valid for all
$0\le a\le 1/2$, to show that whenever $a\le 1/(2k)$
\begin{eqnarray*}
\psi(t)&\ge&
1+t+t^2-\frac{1}{2^{k-1}(1-t)\left(1-\frac{(k-1)t}{2}\right)}\\
&\ge&
1+t^2-\frac{1}{2^{k-1}}(1+t+2t^2)\left(1+\frac{(k-1)t}{2}+2\frac{(k-1)^2t^2}{4}\right)\\
&\ge& 1-\frac{1}{2^{k-1}}+t-\frac{(k+1)t}{2^k}+\frac{t^2}{2}
\enspace,
\end{eqnarray*}
for all $k$ sufficiently large.

The reverse inequality is just as simple
\begin{eqnarray*}
\psi(t)&\le&
1+t+2t^2-\frac{1}{2^{k-1}}(1+t)\left(1+\frac{t}{2}\right)^{k-1}\\
&\le&
1+t+2t^2-\frac{1}{2^{k-1}}(1+t)\left(1+\frac{(k-1)t}{2}\right)\\
&\le& 1-\frac{1}{2^{k-1}}+t-\frac{(k+1)t}{2^k}+2t^2 \enspace.
\end{eqnarray*}
\end{proof}

We are now in position to conclude the proof of Fact
\ref{fact:etae}. Since $\psi$ is increasing and
$\psi({\e_0})=1-\frac{y}{2^{k-1}}$, the inequalities in
(\ref{eq:e}) will be proved once we show that
\begin{eqnarray}\label{eq:need}
\psi\left(\frac{2(1-y)}{2^{k}-k-1}-\frac{16(1-y)^2}{2^{2k}}\right)\le
1-\frac{y}{2^{k-1}}\le \psi\left(\frac{2(1-y)}{2^k-k-1}\right)
\enspace.
\end{eqnarray}
To prove the right-hand inequality in (\ref{eq:need}), set
$t=\frac{2(1-y)}{2^k-k-1}$ and observe that for $k$ large enough,
$t\le 1/(2k)$. Hence, by Lemma \ref{lem:psi},
\begin{eqnarray*}
\psi(t)\ge
1-\frac{1}{2^{k-1}}+t-\frac{(k+1)t}{2^k}=1-\frac{1}{2^{k-1}}+\frac{2(1-y)}{2^k-k-1}-\frac{(k+1)}{2^k}\cdot
\frac{2(1-y)}{2^k-k-1}=1-\frac{y}{2^{k-1}} \enspace.
\end{eqnarray*}

The left-hand inequality in (\ref{eq:need}) is equally simple. In
this case we apply Lemma \ref{lem:psi} with
$t=\frac{2(1-y)}{2^{k}-k-1}-\frac{16(1-y)^2}{2^{2k}}$ and get that
\begin{eqnarray*}
\psi(t)&\le& 1-\frac{1}{2^{k-1}}+t-\frac{(k+1)t}{2^k}+2t^2\\
&\le&
1-\frac{1}{2^{k-1}}+\frac{2(1-y)}{2^k-k-1}-\frac{(k+1)}{2^k}\cdot
\frac{2(1-y)}{2^k-k-1}-\frac{16(1-y)^2}{2^{2k}}+2\cdot
\frac{4(1-y)^2}{(2^k-k-1)^2}\\&\le& 1-\frac{y}{2^{k-1}} \enspace,
\end{eqnarray*}
as long as $k$ is sufficiently large.

\smallskip

To prove the estimate (\ref{eq:eta}) observe that the function
$s\mapsto s(2-s)^{k-1}$ is increasing on $[0,2/k]$. Since we have
shown that for sufficiently large $k$, ${\e_0}\le
\frac{2(1-y)}{2^{k}-k-1}\le \frac{2}{k}$, the lower bound in
(\ref{eq:e}) yields
\begin{eqnarray*}
\eta_0&=&1-{\e_0}(2-{\e_0})^{k-1}\\
&\le&1-2^{k-1}\left(\frac{2(1-y)}{2^{k}-k-1}-\frac{16(1-y)^2}{2^{2k}}\right)
\left(1-\frac{1-y}{2^{k}-k-1}-\frac{8(1-y)^2}{2^{2k}}
\right)^{k-1}\\
&\le& 1-\left(\frac{2^k(1-y)}{2^k-k-1}-\frac{8(1-y)^2}{2^k}\right)
\left(1-\frac{(k-1)(1-y)}{2^{k}-k-1}-\frac{8(k-1)(1-y)}{2^{2k}}
\right)\\
&<& y-\frac{(k+1)(1-y)}{2^k-k-1}+\frac{4k(1-y)^2}{2^k} \enspace,
\end{eqnarray*}
provided $k$ is large enough. The inequality $\eta_0\le y$ is
simpler. By (\ref{eq:identityy}),
$$
\frac{y}{2^k}=\frac{1-{\e_0}(2-{\e_0})^{k-1}}{2(1-{\e_0})(2-{\e_0})^{k-1}}=\frac{\eta_0}{2^k(1-\e_0)(1-\e_0/2)^{k-1}}\ge
\frac{\eta_0}{2^k} \enspace .
$$

\section{Proof of Lemma~\ref{lem:biggerthanhalf}}

Since the function $\alpha\mapsto
\alpha^\alpha(1-\alpha)^{1-\alpha}$ is symmetric around $1/2$, it
suffices to prove that for every $x\in (0,1/2]$,
\begin{eqnarray}\label{eq:notsym}
f\left(\frac12+x,\gamma,\eta\right)>f\left(\frac12-x,\gamma,\eta\right).
\end{eqnarray}
To this end, fix $x\in [-1/2,1/2]$ and $\gamma,\eta>0$. Denote
$\e=1-\gamma^2$. Plugging this notation and $\alpha=1/2+x$ into
(\ref{eq:deff}), we find that the following identity holds
\begin{eqnarray}\label{eq:binom}
\eta^{y/2^{k-1}}2^{2k}\gamma^{2k}f\left(\frac12+x,\gamma,\eta\right)&=&[2x\e^2+(2-\e)^2]^k-2(1-\eta)[(2-\e)+
2x\e]^k+(1-\eta)^2(1+2x)^k\nonumber\\
&=&\sum_{j=0}^k\binom{k}{j}2^jx^{j}[\e^{2j}(2-\e)^{2(k-j)}-2(1-\eta)\e^j(2-\e)^{k-j}+(1-\eta)^2]\nonumber\\
&=&\sum_{j=0}^k\binom{k}{j}2^jx^{j}[\e^j(2-\e)^{k-j}-(1-\eta)]^2.
\end{eqnarray}
This shows that we can write
$f(1/2+x,\gamma,\eta)=\sum_{j=0}^ka_jx^j$ for some $a_j\ge 0$,
such that at most one of the $a_j$'s is zero. Since for every
$x>0$ and odd $j$, $x^j-(-x)^j>0$, (\ref{eq:notsym}) follows.

\section{Proof of Lemma~\ref{lem:decreasing}}

In this section we will use the normalization
(\ref{eq:normalization}). To prove the first assertion of Lemma
\ref{lem:decreasing}, our goal is to show that $g_0'(\alpha) < 0$
for $\frac12 < \alpha \le 1-\frac{3\log k}{k}$. Observe that
\begin{eqnarray}\label{eq:derg0}
g'_0(\alpha)=\frac{f_0(\alpha)^{r-1}\Big\{rf'_0(\alpha)+
f_0(\alpha)\big[\log(1-\alpha)-\log\alpha\big]\Big\}}{\alpha^\alpha(1-\alpha)^{1-\alpha}}
\enspace .
\end{eqnarray}
Differentiating (\ref{eq:f0}) at $x=0$ we find that
$$
f_0'\left(\frac12\right)=2k{\e_0}^2(2-{\e_0})^{2k-2}-4{\e_0}^2(2-{\e_0})^{2k-2}+2{\e_0}^2(2-{\e_0})^{2k-2}=0\enspace.
$$
Since, by (\ref{eq:binom}), $f_0(\alpha)>0$ it is enough to show
that the following function is decreasing on
$\left[\frac12,1-\frac{3\log k}{k}\right]$
$$
\psi(\alpha)=rf_0'(\alpha)+f_0(\alpha)[\log (1-\alpha)-\log\alpha]
\enspace .
$$
Now,
$$
\psi'(\alpha)= rf_0''(\alpha)+f_0'(\alpha)[\log
(1-\alpha)-\log\alpha]-f_0(\alpha)\left(\frac{1}{\alpha}+\frac{1}{1-\alpha}\right)
\enspace .
$$
Since for $1/2< \alpha\le 1$, $\log(1-\alpha)< \log \alpha$ and,
by (\ref{eq:binom}), $f_0'>0$ on $(1/2,1]$, it is thus enough to
prove that
$$
rf_0''(\alpha)\le
f_0(\alpha)\left(\frac{1}{\alpha}+\frac{1}{1-\alpha}\right)
\enspace.
$$
Now, $\frac{1}{\alpha}+\frac{1}{1-\alpha}\ge 4$ and, from
(\ref{eq:f1/2}), we get that for $\alpha\ge 1/2$,
$$
f_0(\alpha)\ge
f_0\left(\frac12\right)=4(1-{\e_0})^2(2-{\e_0})^{2k-2}\ge
(2-{\e_0})^{2k-2} \enspace,
$$
where we also used that, by (\ref{eq:e}), ${\e_0}\le 1/2$ for $k$
large enough. Thus, it suffices to prove
\begin{eqnarray}\label{eq:goal}
rf_0''\left(\frac{1}{2}+x\right)\le 4(2-{\e_0})^{2k-2}.
\end{eqnarray}

Now, using that $x\le\frac12- \frac{3\log k}{k}$, we differentiate
(\ref{eq:f0}) twice to get
\begin{align*}\label{eq:crude}
\lefteqn{f_0''\left(\frac{1}{2}+x\right)} \nonumber\\
&=4k(k-1)\left\{
{\e_0}^4[2x{\e_0}^2+(2-{\e_0})^2]^{k-2}-2{\e_0}^3(2-{\e_0})^{k-1}[2-{\e_0}+2x{\e_0}]^{k-2}+{\e_0}^2(2-{\e_0})^{2k-2}(1+2x)^{k-2}\right\}\nonumber\\
&\le 4k^2\left\{
{\e_0}^4(2-{\e_0})^{2k-4}\left(1+\frac{2x{\e_0}^2}{(2-{\e_0})^2}\right)^{k-2}
+{\e_0}^2(2-{\e_0})^{2k-2}(1+2x)^{k-2}\right\}\nonumber\\
&\le 4k^2\left\{ {\e_0}^4(2-{\e_0})^{2k-4}\left(1+2x\right)^{k-2}
+{\e_0}^2(2-{\e_0})^{2k-2}(1+2x)^{k-2}\right\}\\
&\le 8k^2{\e_0}^2(2-{\e_0})^{2k-2}(1+2x)^k\\ &\le
8k^2{\e_0}^2(2-{\e_0})^{2k-2}\left[2-\frac{6\log k}{k}\right]^k\\
&\le 8k^2{\e_0}^2(2-{\e_0})^{2k-2}2^k\frac{1}{k^3}\\
&\le
\frac{2^{k+3}}{k}\left(\frac{4(1-y)}{2^k}\right)^2(2-{\e_0})^{2k-2},
\end{align*}
where in the last line we used the fact that for $k$ large enough,
(\ref{eq:e}) implies  ${\e_0}\le {4(1-y)}/{2^k}$.

Combining this estimate with (\ref{eq:goal}), we see that we must
show that for sufficiently large $k$
$$
\frac{128\log 2}{1-y+y\log y}\cdot \frac{(1-y)^2}{k}\le 4 \enspace
,
$$
and this is indeed the case by (\ref{eq:stupid}).

\smallskip

It remains to show that $g_0''(1/2)<0$. Denoting
$\zeta(\alpha)=\alpha^{-\alpha}(1-\alpha)^{\alpha-1}$, we see from
(\ref{eq:derg0}) that
$g_0'(\alpha)=f_0(\alpha)^{r-1}\psi(\alpha)\zeta(\alpha)$. Since
$f_0(1/2)=0$, $\zeta'(1/2)=0$, and we have just verified that
$\psi'(1/2)<0$, the required result follows.

\section{Proof of Lemma~\ref{lem:bigalpha}}

Our goal is to show that for any $0<r\le \ttk$, and $1-\frac{3\log
k}{k}< \alpha\le 1$,
\begin{equation}\label{eq:maingoal}
\left[\frac{f(\alpha,\sqrt{\gamma_0},\sqrt{\eta_0})}{f(1/2,\gamma_0,\eta_0)}\right]^r\le
2\alpha^\alpha(1-\alpha)^{1-\alpha} \enspace .
\end{equation}

The following lemma gives an upper bound for the left-hand side of
(\ref{eq:maingoal}).

\begin{lemma}\label{lem:ratio} For all
sufficiently large $k$,
\begin{eqnarray}\label{eq:final}
\frac{f(\alpha,\sqrt{\gamma_0},\sqrt{\eta_0})}{f(1/2,\gamma_0,\eta_0)}\le
y^{y/2^k}\left[1+\frac{2(1-y)-2(1-\sqrt{y})+(1-\sqrt{y})^2\alpha^k}{2^k}+\frac{60k(1-y)^2}{2^{2k}}\right]
\enspace .
\end{eqnarray}
\end{lemma}

\begin{proof}

Denote $\e_1=1-(\sqrt{\gamma_0})^2=1-\sqrt{1-\e_0}$. For
$1-\frac{4\log k}{k}< \alpha\le 1$ write $x=\alpha-1/2$.
Analogously to (\ref{eq:binom}) we have we have the following
identity
\begin{eqnarray}\label{eq:identitye}
&&\!\!\!\!\!\!\!\!\!\!\!\!\eta_0^{y/2^k}2^{2k}\gamma_0^{k}f\left(\frac12+x,\sqrt{\gamma_0},\sqrt{\eta_0}\right)\nonumber\\
&=&
[2x{\e_1}^2+(2-{\e_1})^2]^k-2(1-\sqrt{\eta_0})[(2-{\e_1})+2x{\e_1}]^k+(1-\sqrt{\eta_0})^2(1+2x)^k.
\end{eqnarray}
Our first goal is replace $\eta_0$ in the right-hand side of
(\ref{eq:identitye}) by its upper bound from (\ref{eq:eta}). To
this end consider the function
\begin{eqnarray}\label{eq:rho}
\rho(b)=[2x{\e_1}^2+(2-{\e_1})^2]^k-2b[(2-{\e_1})+2x{\e_1}]^k+b^2(1+2x)^k,
\end{eqnarray}
and observe the right-hand side of (\ref{eq:identitye}) equals
$\rho\left(1-\sqrt{\eta_0}\right)$. So, that it is enough to show
that $\rho$ is decreasing on $[0,1]$. Since $\rho$ is convex and
quadratic, this would follow once we show that $\rho'(1)\le 0$.
This is equivalent to $1+2x\le 2-{\e_1}+2x{\e_1}$, which is true
since $x\le 1/2$. Hence
\begin{eqnarray}\label{eq:upperf}
&&\!\!\!\!\!\!\!\!\!\!\!\!\eta_0^{y/2^k}2^{2k}\gamma_0^{k}f\left(\alpha,\sqrt{\gamma_0},\sqrt{\eta_0}\right)=
\rho(1-\sqrt{\eta_0})
\nonumber\\
&\le& \rho(1-\sqrt{z})=
[2x{\e_1}^2+(2-{\e_1})^2]^k-2(1-\sqrt{z})[(2-{\e_1})+2x{\e_1}]^k+(1-\sqrt{z})^2(1+2x)^k,
\end{eqnarray}
Where $z$ is the upper bound for $\eta_0$ from (\ref{eq:eta}),
i.e.,
\begin{eqnarray}\label{eq:defz}
z=
\min\left\{y,y-\frac{(k+1)(1-y)}{2^k-k-1}+\frac{4k(1-y)^2}{2^k}\right\}
\enspace .
\end{eqnarray}

Hence using $\eta_0\le z$ and the identity (\ref{eq:f1/2}) we
bound the ratio in (\ref{eq:maingoal}) as follows
\begin{align}
\lefteqn{\frac{f(\alpha,\sqrt{\gamma_0},\sqrt{\eta_0})}{f(1/2,\gamma_0,\eta_0)}} \nonumber \\
& =\eta_0^{y/2^k}\cdot\gamma_0^k\cdot
\frac{\eta_0^{y/2^k}2^{2k}\gamma_0^{k}f\left(\alpha,\sqrt{\gamma_0},\sqrt{\eta_0}\right)}
{\eta_0^{y/2^{k-1}}2^{2k}\gamma_0^{2k}f\left(1/2,\gamma_0,\eta_0\right)}\nonumber\\
&=\frac{\eta_0^{y/2^k}(1-\e_0)^{k/2}}{4(1-{\e_0})^2(2-{\e_0})^{2k-2}}\cdot
\eta_0^{y/2^k}2^{2k}\gamma_0^{k}f\left(\alpha,\sqrt{\gamma_0},\sqrt{\eta_0}\right)\nonumber\\
&\le
\frac{z^{y/2^k}(1-\e_0)^{k/2}}{4(1-{\e_0})^2(2-{\e_0})^{2k-2}}\cdot
\left\{[2x{\e_1}^2+(2-{\e_1})^2]^k-2(1-\sqrt{z})[(2-{\e_1})+2x{\e_1}]^k+(1-\sqrt{z})^2(1+2x)^k\right\}\nonumber\\
&= z^{y/2^k}\left[1+\frac{{\e_0}}{2
(1-{\e_0})}\right]^2\left[\frac{\sqrt{1-\e_0}}{(1-\e_0/2)^2}\right]^k\cdot\nonumber\\
&\phantom{\le}
\left\{\left[\frac{x\e_1^2}{4}+\left(1-\frac{\e_1}{2}\right)^2\right]^k-\frac{2(1-\sqrt{z})
[1-{\e_1}(1-\alpha)]^k}{2^k}+\frac{(1-\sqrt{z})^2\alpha^k}{2^k}\right\}
\enspace . \label{eq:estimate}
\end{align}
%

We will bound the various terms in (\ref{eq:estimate}) separately.
First of all, using (\ref{eq:defz}) and the inequality $e^a\le
1+a+a^2$, which is valid for $0\le a\le 1$, we get that
\begin{eqnarray}\label{eq:zpart}
z^{y/2^k}&\le&
y^{y/2^k}\left[1-\frac{(k+1)(1-y)}{y(2^k-k-1)}+\frac{4k(1-y)^2}{y2^k}\right]^{y/2^k}\nonumber\\
&\le&
y^{y/2^k}\exp\left[-\frac{(k+1)(1-y)}{2^k(2^k-k-1)}+\frac{4k(1-y)^2}{2^{2k}}\right]\nonumber\\
&\le&y^{y/2^k}\left[1-\frac{(k+1)(1-y)}{2^k(2^k-k-1)}+\frac{8k(1-y)^2}{2^{2k}}\right],
\end{eqnarray}
as long as $k$ is large enough.

Next, using the inequality $1/(1-a)\le 1+2a$, valid for $0\le a\le
1/2$ we get
\begin{eqnarray}\label{eq:thesquare}
\left[1+\frac{{\e_0}}{2
(1-{\e_0})}\right]^2\le\left[1+\frac{{\e_0}}{2}(1+2{\e_0})\right]^2\le
1+\e_0+5\e_0^2.
\end{eqnarray}
Next, using the inequality $\sqrt{1-x}\le 1-x/2$, the inequality
$1/(1-a)\le 1+a+2a^2$, valid for $0\le a\le 1/2$, and the
inequality $(1+a)^k\le 1+ka+k^2a^2/2$, which is valid for all
$a\le 1/(4k^2)$, we get that since for $k$ large enough $\e_0\le
1/(4k^2)$,
\begin{eqnarray}\label{eq:theroot}
\left[\frac{\sqrt{1-\e_0}}{(1-\e_0/2)^2}\right]^k\le
\frac{1}{(1-\e_0/2)^k}\le
\left(1+\frac{\e_0}{2}+\frac{\e_0^2}{2}\right)^k\le
1+\frac{k\e_0}{2}+\frac{k^2\e_0^2}{8}+k\e_0^2.
\end{eqnarray}
Hence, for $k$ large enough
\begin{eqnarray}\label{eq:theproduct}
\left[1+\frac{{\e_0}}{2
(1-{\e_0})}\right]^2\left[\frac{\sqrt{1-\e_0}}{(1-\e_0/2)^2}\right]^k\le
1+\left(1+\frac{k}{2}\right)\e_0+\frac{k^2\e_0^2}{4}+k\e_0^2.
\end{eqnarray}

Next, using the inequality $x/2\le 1-\sqrt{1-x}\le x/2+x^2$, which
is valid for $0\le x\le 1/2$, we get that
\begin{eqnarray}\label{eq:ineqe1}
\frac{\e_0}{2}\le\e_1=1-\sqrt{1-\e_0}\le \frac{\e_0}{2}+\e_0^2\le
\e_0.
\end{eqnarray}
Observe that since $0\le x\le 1/2$ and $\e_1<1/2$, the function
$\e_1\mapsto \frac{x\e_1^2}{4}+\left(1-\frac{\e_1}{2}\right)^2$ is
decreasing in $\e_1$. Hence, the lower bound in (\ref{eq:ineqe1}),
together with another application of the inequality $(1+a)^k\le
1+ka+k^2a^2/2$, valid for all $a\le 1/(4k^2)$, implies that for
sufficiently large $k$
\begin{eqnarray}\label{eq:firstbracket}
\left[\frac{x\e_1^2}{4}+\left(1-\frac{\e_1}{2}\right)^2\right]^k&\le&
\left[\frac{\e_0^2}{32}+\left(1-\frac{\e_0}{4}\right)^2\right]^k\le
\left[1-\frac{\e_0}{2}+\frac{\e_0^2}{10}\right]^k\le
1-\frac{k\e_0}{2}+\frac{k^2\e_0^2}{8}+k\e_0^2 \enspace .
\end{eqnarray}

The second term in the brackets of (\ref{eq:estimate}) appears
with a minus sign, so we bound it from below, using the fact that
$z\le y$ and $\e_1\le\e_0$.
\begin{eqnarray}\label{eq:secondbracket}
\frac{2(1-\sqrt{z}) [1-{\e_1}(1-\alpha)]^k}{2^k}&\ge&\nonumber
\frac{2(1-\sqrt{z})}{2^k}-\frac{2(1-\sqrt{z})k\e_1(1-\alpha)}{2^k}\\&\ge&
\frac{2(1-\sqrt{z})}{2^k}-\frac{2(1-\sqrt{y})k\e_0}{2^k}\nonumber\\
&\ge&
\frac{2(1-\sqrt{z})}{2^k}-\frac{2(1-y)k}{2^k}\cdot\frac{2(1-y)}{2^{k}-k-1}\nonumber\\
&\ge& \frac{2(1-\sqrt{z})}{2^k}-\frac{8k(1-y)^2}{2^{2k}},
\end{eqnarray}
where we have used the upper bound in (\ref{eq:e}).

Combining (\ref{eq:estimate}), (\ref{eq:zpart}),
(\ref{eq:theproduct}), (\ref{eq:firstbracket}) and
(\ref{eq:secondbracket}) we get

\begin{eqnarray}\label{eq:almost}
&&\!\!\!\!\!\!\!\!\!\!\!\!\!\!\frac{f(\alpha,\sqrt{\gamma_0},\sqrt{\eta_0})}{f(1/2,\gamma_0,\eta_0)}\nonumber\\
&\le&y^{y/2^k}\left[1-\frac{(k+1)(1-y)}{2^k(2^k-k-1)}+\frac{8k(1-y)^2}{2^{2k}}\right]
\left[1+\left(1+\frac{k}{2}\right)\e_0+\frac{k^2\e_0^2}{8}+k\e_0^2\right]\cdot\nonumber\\
&\phantom{\le}&\left\{1-\frac{k\e_0}{2}+\frac{k^2\e_0^2}{8}+k\e_0^2-
\frac{2(1-\sqrt{z})}{2^k}+\frac{8k(1-y)^2}{2^{2k}}+\frac{(1-\sqrt{z})^2\alpha^k}{2^k}\right\}\nonumber\\
&\le&
y^{y/2^k}\left[1-\frac{(k+1)(1-y)}{2^k(2^k-k-1)}+\frac{8k(1-y)^2}{2^{2k}}\right]\cdot\nonumber\\
&\phantom{\le}&\left\{1+\e_0+\frac{(1-\sqrt{z})^2\alpha^k}{2^k}-\frac{2(1-\sqrt{z})}{2^k}
+\frac{8k(1-y)^2}{2^{2k}}+\frac{k(1-\sqrt{z})^2\e_0}{2^k}+\frac{8k^2(1-y)^2\e_0}{2^{2k}}+k^3\e_0^3\right\}\nonumber\\
&\le&y^{y/2^k}\left[1-\frac{(k+1)(1-y)}{2^k(2^k-k-1)}+\frac{8k(1-y)^2}{2^{2k}}\right]\left\{1+\frac{2(1-y)}{2^k-k-1}+\frac{(1-\sqrt{z})^2\alpha^k}{2^k}-\frac{2(1-\sqrt{z})}{2^k}+
\frac{30k(1-y)^2}{2^{2k}}\right\}\nonumber\\
&\le&y^{y/2^k}\left\{1+\frac{2(1-y)}{2^k-k-1}-\frac{(k+1)(1-y)}{2^k(2^k-k-1)}+
\frac{(1-\sqrt{z})^2\alpha^k}{2^k}-\frac{2(1-\sqrt{z})}{2^k}+
\frac{50k(1-y)^2}{2^{2k}}\right\},
\end{eqnarray}
where we have used the upper bound in (\ref{eq:e}),
(\ref{eq:defz}) and the fact that $k$ is large enough.

Now, we claim that for every $\alpha\in [0,1]$,
\begin{eqnarray}\label{eq:trick}
(1-\sqrt{z})^2\alpha^k -2(1-\sqrt{z})\le (1-\sqrt{y})^2\alpha^k
-2(1-\sqrt{y})+z-y.
\end{eqnarray}
Indeed, since by (\ref{eq:defz}), $z\le y$, the left-hand side
minus the right-hand side of (\ref{eq:trick}) is an increasing
function, which vanishes at $1$. Moreover, by (\ref{eq:defz}),
$$
z-y\le -\frac{(k+1)(1-y)}{2^k-k-1}+\frac{4k(1-y)^2}{2^k},
$$
so that (\ref{eq:trick}) becomes
$$ (1-\sqrt{z})^2\alpha^k
-2(1-\sqrt{z})\le (1-\sqrt{y})^2\alpha^k
-2(1-\sqrt{y})-\frac{(k+1)(1-y)}{2^k-k-1}+\frac{4k(1-y)^2}{2^k}.
$$
Plugging this into (\ref{eq:almost}) we get
\begin{eqnarray}\label{eq:almostalmost}
&&\!\!\!\!\!\!\!\!\!\!\!\!\!\!\frac{f(\alpha,\sqrt{\gamma_0},\sqrt{\eta_0})}{f(1/2,\gamma_0,\eta_0)}\le
y^{y/2^k}
\left\{1+\frac{2(1-y)}{2^k-k-1}-\frac{2(k+1)(1-y)}{2^k(2^k-k-1)}+\frac{(1-\sqrt{y})^2\alpha^k}{2^k}-
\frac{2(1-\sqrt{y})}{2^k}+\frac{60k(1-y)^2}{2^{2k}}\right\}\nonumber\\
&=& y^{y/2^k}
\left\{1+\frac{2(1-y)+(1-\sqrt{y})^2\alpha^k-2(1-\sqrt{y})}{2^k}+\frac{60k(1-y)^2}{2^{2k}}\right\}.
\end{eqnarray}

This concludes the proof of Lemma \ref{lem:ratio}.\bigskip
\end{proof}

\mbox{}

 Denote $h(\alpha)=-\alpha\log
\alpha-(1-\alpha)\log(1-\alpha)$. Taking logarithms of
(\ref{eq:maingoal}), and using (\ref{eq:final}) and the inequality
$\log (1+x)\le x$, we see that our goal is reduced to showing that
for all $r\le \ttk$,
\begin{eqnarray}\label{eq:newgoal}
\frac{r}{2^k}&&\!\!\!\!\!\!\!\!\!\!\!\left[y\log
y+2(1-y)-2(1-\sqrt{y})+(1-\sqrt{y})^2\alpha^k+\frac{60k(1-y)^2}{2^{k}}\right]
\le \log 2-h(\alpha) \enspace .
\end{eqnarray}

For simplicity denote:
\begin{eqnarray}\label{eq:abc}
\left\{\begin{array}{ll} A=(1-\sqrt{y})^2 \\
B=y\log
y+2(1-y)-2(1-\sqrt{y})+\frac{60k(1-y)^2}{2^{k}}\end{array}\right.
\end{eqnarray}
With this notation (\ref{eq:newgoal}) becomes
$$
\frac{r}{2^k}\le \frac{\log2-h(\alpha)}{A\alpha^k+B}\equiv
M(\alpha) \enspace ,
$$
and this should hold for all $\alpha \ge 1-\frac{3\log k}{k}$. We
are therefore interested in the minimal value of $M$ on the
interval $\left[1-\frac{3\log k}{k},1\right]$. The derivative of
$M$ is
\begin{eqnarray}\label{eq:derM}
M'(\alpha)=\frac{(A\alpha^k+B)\cdot[\log
\alpha-\log(1-\alpha)]-kA\alpha^{k-1}[\log2-h(\alpha)]}{(A\alpha^k+B)^2}.
\end{eqnarray}
In particular, $M'(1)=\infty$, so that the minimum of $M$ cannot
occur at $\alpha=1$. We rule out the possibility of the minimum
being at $1-\frac{3\log k}{k}$ in the following claim.
\begin{claim}\label{claim:notleft} If $k$ is large enough
then for every $1-\frac{3\log k}{k}\le \alpha\le 2^{-1/k}$,
$M\left(\alpha\right)>M(1)$.
\end{claim}
\begin{proof} Observe that for every $\beta\in [0,1]$,
\begin{eqnarray}\label{eq:ineqh}
h(\beta)&=&\beta\log(1/\beta)-(1-\beta)\log(1-\beta)\nonumber\\&\le&
\beta\left(\frac{1}{\beta}-1\right)-(1-\beta)\log(1-\beta)=1-\beta-(1-\beta)\log(1-\beta)
\enspace .
\end{eqnarray}
Hence, since $\alpha\ge 1-\frac{3\log k}{k}$,
\begin{eqnarray}\label{eq:hleft}
h(\alpha)\le h\left(1-\frac{3\log k}{k}\right)\le \frac{3\log
k}{k}+\frac{3\log k}{k}\log\left(\frac{k}{3\log k}\right)\le
\frac{4(\log k)^2}{k} \enspace .
\end{eqnarray}
Using the fact that $\alpha\le 2^{-1/k}$, it follows that
\begin{eqnarray*}
M\left(\alpha\right)\ge \frac{\log 2-\frac{5(\log
k)^2}{k}}{A\left(2^{-1/k}\right)^k+B}\ge \frac{\log
2\left(1-\frac{10(\log k)^2}{k}\right)}{\frac{A}{2}+B} \enspace .
\end{eqnarray*}
On the other hand, $M(1)=\log 2/(A+B)$, so that it is enough to
show that
$$
1-\frac{10(\log
k)^2}{k}\ge\frac{\frac{A}{2}+B}{A+B}=1-\frac{\frac{1}{2}A}{A+B}
\enspace ,
$$
which is equivalent to
\begin{eqnarray}\label{eq:stupidgoal}
\frac{A}{A+B}\ge \frac{20(\log k)^2}{k} \enspace .
\end{eqnarray}
Observe that since $1-\sqrt{y}\ge (1-y)/2$, $A\ge (1-y)^2/4$. On
the other hand, using (\ref{eq:stupid}) we get that for
sufficiently large $k$,
\begin{eqnarray*}
A+B=1-y+y\log y+\frac{60k(1-y)^2}{2^{k}} \le
(1-y)^2+\frac{60k(1-y)^2}{2^k}\le 2(1-y)^2 \enspace ,
\end{eqnarray*}
It follows that the left-hand side in (\ref{eq:stupidgoal}) is at
least $1/8$, so that (\ref{eq:stupidgoal}) provided $k$ large
enough.
\end{proof}

By Claim \ref{claim:notleft} it remains to bound $M(\alpha)$ from
below when $\alpha>2^{-1/k}$ and $M'(\alpha)=0$. In this case, by
(\ref{eq:derM}),
\begin{eqnarray}\label{eq:condition}
-\log(1-\alpha)= -\log\alpha+
\frac{kA\alpha^{k-1}}{A\alpha^k+B}[\log 2-h(\alpha)] \enspace .
\end{eqnarray}
From the lower bound $\alpha\ge 2^{-1/k}$ and (\ref{eq:ineqh}) it
follows that $h(\alpha)\le \frac{4\log k}{k}$. Hence
(\ref{eq:condition}), together with our assumption that $k$ is
large, implies
\begin{eqnarray}\label{eq:alphaisbig}
-\log(1-\alpha)>\frac{kA/2}{A+B}\left[\log2-\frac{4\log
k}{k}\right]\ge \frac{k}{5}\cdot\frac{A}{A+B} \enspace .
\end{eqnarray}
As we have seen in the proof of Claim \ref{claim:notleft}, $A\ge
(1-y)^2/4$ and $A+B\le 2(1-y)^2$. Plugging these inequalities into
(\ref{eq:alphaisbig}), we get that $-\log(1-\alpha)>k/40$, i.e.
$\alpha\ge 1-e^{-k/40}$. Plugging this into (\ref{eq:condition})
once more, we get that
\begin{eqnarray*}
-\log(1-\alpha)\ge \frac{kA[1-2k/(e^{k/40})]}{A+B}[\log
2-2/(e^{k/40})]\ge \frac{kA\log
2}{A+B}\left(1-\frac{6}{e^{k/40}}\right) \enspace .
\end{eqnarray*}
Finally, we have shown that
\begin{eqnarray}\label{eq:maximizer}\alpha\ge 1-\exp\left[-\frac{kA\log
2}{A+B}\left(1-\frac{6}{e^{k/40}}\right)\right] \enspace .
\end{eqnarray}

We are now ready to bound $M(\alpha)$ from below. We start by
recalling that
\begin{eqnarray*}
A+B= 1-y+y\log y+\frac{60k(1-y)^2}{2^{k}} \equiv 1-y+y\log y+P
\enspace .
\end{eqnarray*}
Using the inequality $1/(1+x)\ge 1-x$, we get
\begin{eqnarray}\label{eq:A+C}
\frac{1}{A+B}\ge \frac{1}{1-y+y\log y}\left[1-\frac{P}{1-y+y\log
y}\right] \ge\frac{1}{1-y+y\log y}\left[1-\frac{120k}{2^k}\right],
\end{eqnarray}
where the last inequality used (\ref{eq:stupid}). Of course, we
also know that $A+B\ge 1-y+y\log y$.

Now, using (\ref{eq:ineqh}), (\ref{eq:maximizer}), (\ref{eq:A+C}),
and the fact that $(1-\sqrt{y})^2/(1-y+y\log y)\le 1$, we get
\begin{eqnarray*}
h(\alpha)&\le& \frac{2kA\log
2}{A+B}\left(1-\frac{6}{e^{k/40}}\right)\exp\left[-\frac{kA\log
2}{A+B}\left(1-\frac{6}{e^{k/40}}\right)\right]\\
&\le& 2k\log2\exp\left[-\frac{k(1-\sqrt{y})^2\log 2}{1-y+y\log
y}\left(1-\frac{120k}{2^k}\right)\left(1-\frac{6}{e^{k/40}}\right)\right]\\
&\le&2k\log2\exp\left[-\frac{k(1-\sqrt{y})^2\log 2}{1-y+y\log
y}\left(1-\frac{120k}{2^k}-\frac{6}{e^{k/40}}\right)\right]\\
&\le& 10k\log 2 \cdot 2^{-k\varphi(y)},
\end{eqnarray*}
where  as in Proposition \ref{prop:main},
$\varphi(y)=\frac{(1-\sqrt{y})^2}{1-y+y\log y}$.

\smallskip
So, using (\ref{eq:A+C}), we get
\begin{eqnarray*}
M(\alpha)&=&\frac{\log2-h(\alpha)}{A\alpha^k+B}\\
&\ge& \frac{\log 2}{A+B}(1-10k2^{-k\varphi(y)})\\
&\ge&\frac{\log 2}{1-y+y\log
y}\left[1-\frac{120k}{2^k}\right](1-10k2^{-k\varphi(y)})\\
&\ge&\frac{\log 2}{1-y+y\log
y}\left[1-\frac{120k}{2^k}-10k2^{-k\varphi(y)}\right]\\
&\ge& \frac{\log 2}{1-y+y\log
y}\left[1-20k2^{-k\varphi(y)}\right],
\end{eqnarray*}
where we have used the fact that $\varphi(y)\ge 1/2$ and that $k$
is sufficiently large.

\medskip

This concludes the proof of Lemma \ref{lem:bigalpha}.

\end{document}